\documentclass[hidelinks,onefignum,onetabnum]{siamart251216}

\usepackage{graphicx}
\ifpdf
  \DeclareGraphicsExtensions{.eps,.pdf,.png,.jpg}
\else
  \DeclareGraphicsExtensions{.eps}
\fi

\usepackage{amsmath}
\usepackage{amssymb}
\usepackage[margin=1in]{geometry}
\usepackage{multirow}
\usepackage{siunitx}
\usepackage{comment}
\usepackage{tikz}
\usepackage{braket}
\usepackage{float}

\usetikzlibrary{decorations.pathmorphing, decorations.pathreplacing, arrows.meta,positioning,fit}
\usetikzlibrary{decorations.markings}

\usepackage{amsfonts}
\usetikzlibrary{positioning}

\usepackage{algorithm}
\usepackage{algpseudocode}

\usepackage{enumitem}
\usepackage{physics}
\usepackage{mathtools}
\usepackage[ntheorem]{empheq}
\usepackage{titlesec}
\usepackage{subcaption}

\newcommand\bm[1]{\begin{bmatrix}#1\end{bmatrix}} 
\newcommand{\Nnode}{N_q}

\newsiamremark{remark}{Remark}
\newsiamremark{hypothesis}{Hypothesis}
\crefname{hypothesis}{Hypothesis}{Hypotheses}
\newsiamthm{claim}{Claim}
\newsiamremark{fact}{Fact}
\crefname{fact}{Fact}{Facts}

\headers{Gregory NPI Schemes for the Lindblad Equation}{  J. Hu, D. Appel\"{o}, Y. Cheng}

\title{Gregory Nested Picard Iteration Schemes for Open Quantum Systems Governed by the Lindblad Equation\thanks{Submitted to the SISC Special Section on Quantum Computing on \today.}
\funding{
DA is supported by the U.S. Department of Energy, Office of Science, Advanced Scientific Computing Research (ASCR), under Award Number DE-SC0025424. This material is based upon work supported, in part, by the National Science Foundation under Grant No. DMS-2436319 and Virginia Tech. 
YC is supported by DOE grant DE-SC0023164, AFOSR grant FA9550-25-1-0154 and Virginia Tech. This material is based upon work supported by the National Science Foundation under Grant No.
DMS-2424139 while the second and third authors were in residence at the Simons Laufer Mathematical Sciences Institute in Berkeley, California, during the Fall 2025 semester.
}}

\author{
    Jiuhua Hu\thanks{Department of Mathematics, Virginia Tech, Blacksburg, VA 24060 (\email{jiuhuahu115@gmail.com, appelo@vt.edu, yingda@vt.edu}).}
    \and
     Daniel Appel\"{o}\footnotemark[2]    
    \and 
    Yingda Cheng\footnotemark[2]
}

\usepackage{amsopn}

\ifpdf
\hypersetup{
  pdftitle={HO LR SP Lindblad Gregory},
  pdfauthor={J. Hu, D. Appel\"{o}, Y. Cheng}
}
\fi

\begin{document}

\maketitle

%
%
\begin{abstract}
Numerical simulation of quantum computing hardware and open quantum systems governed by the Lindblad equation is challenging due to the high dimensionality of the density matrix and the need to preserve fundamental physical properties. In our previous work \cite{hu2025arbitrary}, we developed an arbitrary-order, low-rank, completely positive and trace preserving (CPTP) method for the Lindblad equation with time-dependent Hamiltonians by nested Picard iteration (NPI). In this work, we  develop Gregory NPI schemes which are CPTP schemes constructed by Gregory-type quadrature on equispaced nodes. The methods, which are of order up to nine, substantially reduce the computational cost compared to NPI schemes with Gaussian quadrature rules in \cite{hu2025arbitrary}, while retaining high-order accuracy and structure preservation. We analyze  the stability of the resulting scheme  for a physics-based test equation. Numerical experiments verify the convergence of the method,  and demonstrate the effectiveness of the low-rank approximation. We study the performance of the CNOT gate constructed in \cite{lee2026high} for both closed and open quantum systems. 
\end{abstract}

%
%


%
%
\section{Introduction}

Open quantum systems refer to quantum systems that interact with their surrounding environment and therefore cannot be regarded as isolated. Such systems govern the hardware used in quantum computers.  In this setting, the dynamics is governed not only by the Hamiltonian but also by dissipative and decoherence effects induced by the environment. A standard model for such dynamics is the Lindblad equation, which describes the evolution of the density matrix $\rho(t)$. 

The Lindblad equation plays a central role in quantum information processing, quantum control, and the modeling of realistic quantum devices, where dissipation and decoherence are unavoidable. Accurate and efficient numerical simulation of the density matrix is therefore of fundamental importance. However, this task is challenging for at least two reasons. First, the state space grows rapidly with the number of subsystems: for example, for a $Q$-qubit system, the density matrix has dimension $2^Q\times 2^Q$. Second, a  numerical method should preserve the key structural  properties of the density matrix, such as Hermiticity, positive semidefiniteness, and unit trace. These requirements make the construction of efficient high-order numerical integrators challenging. 

Existing methods that are designed to preserve the completely positive (CP) property often do so through constructions motivated by Choi's theorem; see, for example, \cite{appelo2024kraus,cao2025structure,choi1975completely,hu2025arbitrary,kraus1983states,robin2025unconditionally}. In particular, a numerical update
\[
\rho^n \mapsto \rho^{n+1}=\mathcal{G}(\rho^n),
\qquad \rho^n \approx \rho(t^n),
\]
is completely positive if the linear map $\mathcal{G}$ admits a Kraus representation of the form
\[
\mathcal{G}(\rho)=\sum_l G_l \rho G_l^{\dagger}.
\]
The unit-trace condition can then be enforced, when necessary, by the renormalization
\[
\rho^{n+1}\leftarrow \frac{\rho^{n+1}}{\operatorname{Tr}(\rho^{n+1})}.
\]

High-order methods for the Lindblad equation may be constructed either stochastically or deterministically. One stochastic approach is based on trajectory methods, often referred to as high-order unravelings. Such methods evolve random pure states, or random density matrices, associated with individual quantum trajectories determined by measurement records or jump events. The Lindblad solution is recovered only after averaging over these trajectories; see, for example, \cite{steinbach1995high,CPTPhighOrderUnraveling2024}.
 Consequently, in addition to time discretization error, they typically incur sampling error. An alternative is provided by deterministic high-order integrators that evolve the density matrix directly. Among them, an important class is based on Duhamel's principle; see, for example, \cite{appelo2024kraus,cao2025structure,LiWang2023HigherOrder,hu2025arbitrary}.

In our previous work \cite{hu2025arbitrary}, we proposed high-order low-rank CPTP schemes based on nested Picard iteration (NPI) with Gaussian quadrature rules. The use of Gaussian quadrature was motivated in part by the positivity of its weights, which is crucial for preserving complete positivity; see also \cite{LiWang2023HigherOrder} which uses Gaussian quadrature for high-dimensional integrals arising in the time-ordered methods of \cite{cao2025structure,LiWang2023HigherOrder}. 
While structure-preserving, the schemes in \cite{hu2025arbitrary} are costly in particular in high order cases due to the need to repeatedly recalculate the density matrix at intermediate stages.  To compute a solution of order $p$, the Gaussian-quadrature-based NPI scheme requires $O(\frac{p!}{2^p})$ flow-operator evaluations.

In this paper, we propose the Gregory NPI schemes to  improve the efficiency of   \cite{hu2025arbitrary} by replacing Gaussian quadrature with Gregory quadrature on equidistant nodes. The resulting schemes remain high-order (up to ninth order), low-rank, and structure-preserving, while being significantly more efficient. In particular, to compute a solution of order $p$,  the scheme proposed here requires only $\Nnode^{(p)}+1$ flow-operator evaluations, where $\Nnode^{(p)}+1$ is the number of Gregory quadrature nodes, which is taken to be $2p-2$.  Therefore, the scheme only requires $O(p)$ flow operator evaluations, which  demonstrates the large reduction in computational cost   particularly for high-order cases.  

The remainder of this paper is organized as follows. In Section~\ref{sec:preliminary}, we review the Lindblad equation and Gregory quadrature. In Section~\ref{sec:Gaussian_NPI_Schemes}, we review the Gaussian-quadrature-based NPI schemes. In Section~\ref{sec:Gregory_NPI_Schemes}, we formulate the new and more efficient method based on Gregory quadrature. In Section~\ref{sec: stability}, we study the stability of the Gregory NPI schemes without trace renormalization for a single-qubit model problem. In Section~\ref{sec:truncation}, we describe the low-rank truncation. In Section~\ref{sec: numerical experiments}, we present numerical results for representative quantum systems. These include a two-qubit problem with a known exact solution, a controlled qudit-resonator system for examining convergence, low-rank truncation, and the effect of trace renormalization, and a CNOT-gate problem involving two qudits and a resonator. We conclude the paper in Section~\ref{sec:conclusion}.

\section{Preliminaries}\label{sec:preliminary}
In this section, we introduce the Lindblad equation and Gregory quadrature. 

\subsection{Governing Equations}
In this work, the evolution of the density matrix $\rho(t)\in\mathbb{C}^{d\times d}$ is governed by the Lindblad equation
\begin{align}\label{eqn_master}
    \frac{d\rho(t)}{dt}
    &= -i\bigl(H(t)\rho(t)-\rho(t)H(t)\bigr)+\mathcal{L}(\rho(t)),\\
    \rho(0) &= \rho_0. \nonumber
\end{align}
Here $\rho_0$ is the initial density matrix, and
\[
H(t)=H_d+H_c(t),
\]
is the system Hamiltonian, where $H_d$ is the drift Hamiltonian and $H_c(t)$ is a time-dependent control Hamiltonian.

The dissipative part of the dynamics is described by the Lindbladian
\begin{equation}\label{eq:Lindbladian_general}
\mathcal{L}(\rho)
=
\sum_{\alpha}\mathcal{L}_{\alpha}\rho \mathcal{L}_{\alpha}^{\dagger}
-\frac{1}{2}\sum_{\alpha}
\left(
\mathcal{L}_{\alpha}^{\dagger}\mathcal{L}_{\alpha}\rho
+
\rho\,\mathcal{L}_{\alpha}^{\dagger}\mathcal{L}_{\alpha}
\right),
\end{equation}
where $\{\mathcal{L}_{\alpha}\}_{\alpha}$ is a family of jump operators and $\dagger$ denotes the adjoint.

Our main interest is methods for simulating quantum computers and thus we consider a transmon-based model for a composite superconducting quantum device; see \cite{devoret2013superconducting,PhysRevLett.101.080502,PhysRevA.76.042319,gunther2021quantum}. We adopt the notation used in \cite{gunther2021quantum}. For a system consisting of $Q$ subsystems, where subsystem $k$ has $n_k$ energy levels, the drift Hamiltonian is given by
\[
H_d
=
\sum_{k=0}^{Q-1}
\left(
\omega_k a_k^\dagger a_k
-\frac{\xi_k}{2}\,a_k^\dagger a_k^\dagger a_k a_k
+\sum_{l>k}
\Big(
J_{kl}(a_k^\dagger a_l+a_k a_l^\dagger)
-\xi_{kl} a_k^\dagger a_k a_l^\dagger a_l
\Big)
\right),
\]
where $\omega_k\ge0$ is the $0\to1$ transition frequency of subsystem $k$, $\xi_k\ge0$ is the self-Kerr coefficient, and $J_{kl}\ge0$ and $\xi_{kl}\ge0$ denote the coupling and cross-Kerr coefficients between subsystems $k$ and $l$, respectively. The lowering operator on subsystem $k$ is defined by
\[
a_k
=
I_{n_0}\otimes\cdots\otimes I_{n_{k-1}}
\otimes a^{(n_k)}
\otimes I_{n_{k+1}}\otimes\cdots\otimes I_{n_{Q-1}},
\]
with
\[
a^{(n_k)}
=
\bm{0&1 & && \\
&0&\sqrt{2} && \\
&&\ddots&\ddots& \\
&&&0&\sqrt{n_k-1}\\
&&&&0}.
\]
The control Hamiltonian is given by
\[
H_c(t)=\sum_{k=0}^{Q-1} f^k(t)(a_k+a_k^\dagger),
\]
where $f^k(t)$ is a real-valued control function. For the transmon model considered here, the jump operators in \eqref{eq:Lindbladian_general} are chosen to model relaxation and pure dephasing. More precisely,
\[
\mathcal{L}_{1k}=\frac{1}{\sqrt{T_1^k}}\,a_k,
\qquad
\mathcal{L}_{2k}=\frac{1}{\sqrt{T_2^k}}\,a_k^\dagger a_k,
\]
where $T_1^k>0$ and $T_2^k>0$ denote the characteristic relaxation and pure dephasing times for subsystem $k$, respectively. Accordingly, the Lindbladian takes the form
\[
\mathcal{L}(\rho)
=
\sum_{k=0}^{Q-1}\sum_{l=1}^2
\mathcal{L}_{lk}\rho \mathcal{L}_{lk}^\dagger
-\frac12
\left(
\mathcal{L}_{lk}^\dagger \mathcal{L}_{lk}\rho
+
\rho \mathcal{L}_{lk}^\dagger \mathcal{L}_{lk}
\right).
\]
The operators $\mathcal{L}_{1k}$ and $\mathcal{L}_{2k}$ represent energy dissipation and dephasing.

\subsection{Classical Gregory Quadrature}
In this paper, we will consider the approximation of the solution to the initial value problem \eqref{eqn_master} from time 0 to time $T$ on an equidistant grid 
\begin{equation}\label{eq:time_grid}
t^n = n \Delta t, \ \ n = 0,\ldots,N_T, \ \ \Delta t = T/N_T.
\end{equation}
In our scheme, we introduce  $\Delta T := N_q^{(p)} \Delta t$ as the ``big" interval where the numerical quadrature will be performed. Here $N_q^{(p)}+1$ is the number of Gregory nodes. In order to achieve $p$-th order accuracy $N_q^{(p)}$ must be at least $2p-3$. This is the value we always choose for $N_q^{(p)}$ in this paper.

We now review the basics of classical Gregory quadrature rule. For a scheme of order $p \geq 2,$ we use $N_q^{(p)} + 1$ nodes to compute an integral $\int_{t^n}^{t^{n+N_q^{(p)}}} g(t) dt$ which is an integral of a smooth function $g(t)$ over the time interval $[t^n,t^{n+N_q^{(p)}}]$ of length $\Delta T$. 
We let $\tau^j = t^n + j\Delta t$ with $j = 0,\ldots,N_q^{(p)}$ and let
\[
g^j:=g(\tau^j),\qquad j=0,1,\dots,N_q^{(p)}.
\]
The classical Gregory rule may be viewed as an endpoint-corrected composite trapezoidal rule, only the weights near the two endpoints are modified, while the interior weights are equal to one. The general definition of quadrature weights is provided in \cite{fornberg2021improving}. In this paper, we use  
\begin{equation}\label{eq:Gregory-weighted}
Q_{N_q^{(p)}}^{\mathrm{G}}[g]
:=
\Delta t \sum_{j=0}^{N_q^{(p)}} \omega_j^{(p)} g^j,
\end{equation}
where the weights $\omega_j^{(p)}, j=0, \ldots p-2$   are provided in Table~\ref{table:gregory}. The weights for $j=p-1, \ldots, 2p-3$ are the same but with the order reversed. Note that with the choice we make here, $N_q^{(p)} := 2p-3$, there are no interior weights.  For the classical Gregory weights, if $p \ge 10,$ there will be negative weights associated with the Gregory quadrature, which will break the desired CP property. Therefore, in this paper, we only consider order below or equal to 9 (we note that even higher-order methods can straightforwardly obtained by using the weights proposed by Fornberg \cite{fornberg2021improving}). 
\begin{table}[] 
\renewcommand{\arraystretch}{1.5}  
\centering 
\caption{Gregory weights associated with nodes near the left endpoint. Note that all of the weights are positive. \label{table:gregory}}
\begin{tabular}{c|l}
\hline
Order $p$  & Weights $\omega_j^{(p)}, j=0, \ldots, p-2$  \\
\hline
2 & \( \frac{1}{2} \) \\
\hline
3 & \( \frac{5}{12}, \, \frac{13}{12} \) \\
\hline
4 & \( \frac{3}{8}, \, \frac{7}{6}, \, \frac{23}{24} \) \\
\hline
5 & \( \frac{251}{720}, \, \frac{299}{240}, \, \frac{211}{240}, \, \frac{739}{720} \) \\
\hline
6 & \( \frac{95}{288}, \, \frac{317}{240}, \, \frac{23}{30}, \, \frac{793}{720}, \, \frac{157}{160} \) \\
\hline
7 & \( \frac{19087}{60480}, \, \frac{84199}{60480}, \, \frac{18869}{30240}, \, \frac{37621}{30240}, \, \frac{55031}{60480}, \, \frac{61343}{60480} \) \\
\hline
8 & \( \frac{5257}{17280}, \, \frac{22081}{15120}, \, \frac{54851}{120960}, \, \frac{103}{70},\,\frac{89437}{120960}, \, \frac{16367}{15120}, \, \frac{23917}{24192} \) \\
\hline
9 & \( \frac{1070017}{3628800}, \, \frac{5537111}{3628800}, \, \frac{103613}{403200}, \, \frac{261115}{145152}, \, \frac{298951}{725760}, \, \frac{515677}{403200}, \, \frac{3349879}{3628800}, \, \frac{3662753}{3628800} \) \\
\hline
\end{tabular}
\end{table}

\section{CPTP Schemes via NPI}\label{sec:Gaussian_NPI_Schemes}
In this section, we  review the numerical framework developed in \cite{hu2025arbitrary} for constructing high-order completely positive trace-preserving  schemes based on NPI. 

Following \cite{steinbach1995high}, the Lindblad equation can be reformulated in the form
\begin{equation}\label{eqn:Lindblad_in_J_term}
\cfrac{d\rho(t)}{dt}
= J(t)\rho(t) + \rho(t)J^\dagger(t)
+ \sum_\alpha \mathcal L_\alpha \rho(t)\mathcal L_\alpha^\dagger ,
\end{equation}
where
\[
J(t) = -i H_{\mathrm{eff}}, 
\qquad 
H_{\mathrm{eff}} = H(t) + \cfrac{1}{2i}\sum_\alpha \mathcal L_\alpha^\dagger \mathcal L_\alpha .
\]
For convenience we introduce the notation 
\[
\mathcal{ L}_L\rho(t)
=
\sum_\alpha
\mathcal L_\alpha \rho(t)\mathcal L_\alpha^\dagger .
\]

To illustrate the main idea from \cite{hu2025arbitrary}, we first consider the simplified equation without \(\mathcal{L}_{L}(\rho(t))\) (which is already on Kraus form) 
\begin{equation}
\frac{d\rho(t)}{dt}
=
J(t)\rho(t)+\rho(t)J^\dagger(t), \ \ \rho(t_0) = \rho_0.
\end{equation}
Since the right-hand side is not expressed in Kraus form, directly applying
a standard ODE solver (for instance an explicit Runge--Kutta method) does not
guarantee preservation of the CPTP structure. However, the density matrix
$\rho(t)$ is positive semidefinite and therefore admits a factorization
\[
\rho(t) = V(t)V^\dagger(t).
\]
Substituting this representation into the above equation shows that solving
for $\rho(t)$ is equivalent to solving the matrix differential equation
\begin{equation}
\frac{dV(t)}{dt} = J(t)V(t), 
\quad
V(t_0)=V_0,\label{eq:J(t)}
\end{equation}
where the initial factor $V_0$ satisfies
\[
\rho_0 = V_0 V_0^\dagger .
\]
Now, the solution of \eqref{eq:J(t)} can be written in terms of the flow operator
$U(t,s)$ associated with the equation, namely
\[
V(t) = U(t,s)V(s).
\]
Consequently the density matrix evolves according to
\begin{equation} \label{eq:dflow}
\rho(t) = U(t,s)\rho(s)U^\dagger(t,s).
\end{equation}

Equation \eqref{eq:dflow} in combination with treating the  Kraus term in the full Lindblad equation as a forcing makes it natural to apply Duhamel's principle. Then the solution of the complete equation \eqref{eqn:Lindblad_in_J_term} can be expressed as
\begin{equation} \label{eq:Duhamel_form}
\rho(t)
=
U(t,\tau)\rho(\tau)U^\dagger(t,\tau)
+
\int_\tau^t
U(t,s)\,\mathcal{L}_L\rho(s)\,U^\dagger(t,s)\,ds.
\end{equation}

\subsection{Nested Picard Iteration using Gaussian Quadrature}
\begin{figure}[]
\centering
\scalebox{0.7}{%
\begin{tikzpicture}[>=Stealth, thick]

\tikzset{
  box/.style={
    draw,
    rounded corners=2pt,
    minimum width=1.8cm,
    minimum height=1.2cm,
    inner sep=3pt
  },
  bigbox/.style={
    draw,
    rounded corners=2pt,
    minimum width=1.8cm,
    minimum height=1.2cm,
    inner sep=6pt
  },
  lab/.style={font=\small, fill=white, inner sep=1pt}
}

\node[bigbox] (B) at (3.5,0) {};
\node at (3.5,-0.05) {$\rho^n$};

\node[box] (p1) at (-2,2.) {$\rho^{n+\frac{1}{3},1}$};
\node[box] (p2) at (3.5,2.5) {$\rho^{n+\frac{2}{3},2}$};
\node[box] (p3) at (8.5,3.5) {$\rho^{n+1,3}$};

\draw[->] ([xshift=-0.5cm]B.north) -- (p1.south)
    node[midway,left,lab] {$U^{(1)}(t^{n+\frac{1}{3}},t^n)$};

\draw[->] (B.north) -- (p2.south)
    node[midway,right,lab] {$U^{(2)}(t^{n+\frac{2}{3}},t^n)$};

\draw[->] ([xshift=0.5cm]B.north) -- (p3.south)
    node[midway,right,lab,align=center] {
    {$U^{(3)}(t^{n+1},t^n)$}\\[4pt]
    {$U^{(2)}(t^{n+1},t^n)$}
    };
\draw[->] (p1.east) -- (p2.west)
    node[midway,above,lab] {$U^{(2)}(t^{n+\frac{2}{3}},t^{n+\frac{1}{3}})$};

\draw[->] (p2.east) -- (p3.west)
    node[midway,above,lab] {$U^{(3)}(t^{n+1},t^{n+\frac{2}{3}})$};

\end{tikzpicture}%
}

\caption{Illustration of the third-order Gaussian-quadrature-based NPI scheme. To advance the solution from $\rho^n$ to $\rho^{n+1}$, two intermediate nodes are required.}
\label{fig:limitation_example_order_3}

\end{figure}
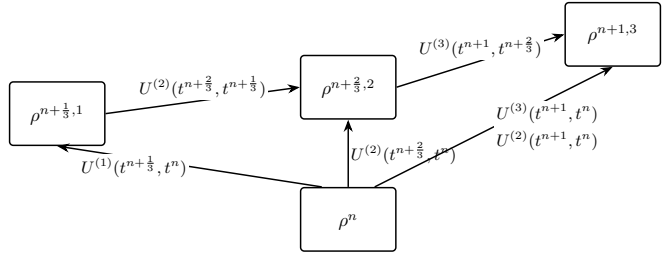
In \cite{hu2025arbitrary} we   discretize \eqref{eq:Duhamel_form}  using NPI. The basic idea is to gain one order of accuracy per iteration. 
Precisely, using the time grid \eqref{eq:time_grid} with timestep $\Delta t$ a solution of $p$-th accuracy can be obtained by the following iteration 
\begin{multline}
\label{eq:NPI_general} 
  \rho^{n+1,k+1} = U^{(k+1)}(t^{n+1},t^n) \rho^n U^{(k+1)}(t^{n+1},t^n) ^\dagger  \\
 +\Delta t \sum_{j=0}^{N_q^{(k+1)}} \omega_j^{(k+1)} U^{(k)}(t^{n+1},t^{n+c^{(k+1)}_j}) \mathcal{L}_L \rho^{n+c^{(k+1)}_j,k} U^{(k)}(t^{n+1},t^{n+c^{(k+1)}_j}) ^\dagger,  
\end{multline}
with $k =0,1,\cdots, p-1$. The iteration is started by $\rho^{n+1,0} = \rho^n$. Here, \(U^{(k)}(t,s)\) represents a $k$-th order approximation of the flow operator $U(t,s)$, and  $\omega_j^{(k+1)}$ and  $c_j^{(k+1)}$ are the weights and nodes of the $(k+1)$-th order Gaussian quadrature rule. Note here we have flexibility to choose various types of Gaussian quadrature (Gaussian, Gauss-Radau, Gauss-Lobatto) which will have a slight effect on  the efficiency of the method.  

As an example, below we show the details of a third order scheme, where the quadrature rule is taken as Gauss-Radau when $k=3$ and midpoint rule for $k=2.$ Figure~\ref{fig:limitation_example_order_3} illustrates the details of the computations involved. The computation of $\rho^{n+1,3}$ requires evaluating intermediate solutions $\rho^{n+\frac{1}{3},1}$ and $\rho^{n+\frac{2}{3},2}$ at the new nodes $t^{n+\frac{1}{3}}$ and $t^{n+\frac{2}{3}}$. This process involves six flow operators of different orders acting over different subintervals. 

We can give an estimate on how the cost scales with order of accuracy. For Gaussian quadrature the number of quadrature nodes, $N_q^{(p)}+1$, to obtain order $p$ is approximately $p/2$ (accounting for the option to include endpoints or not). Let $M(p)$ be the number of Gaussian nodes and let  $F(p)$ be the number of flow-operator evaluations with the sum taken over iteration $1$ to $p$ (i.e. $k = 0,\ldots,p-1$ in \eqref{eq:NPI_general}). Then $M(p)$ and $F(p)$ can be computed by the recurrence relations
\begin{eqnarray*}
M(p+1) &=& \left(M(p)+1\right)\left(N_q^{(p+1)}+1\right) \approx \left(M(p)+1\right) \frac{p+1}{2} ,\\
F(p+1) &=& F(p)\bigl(N_q^{(p+1)}+1\bigr)+N_q^{(p+1)}+2 \approx F(p)\frac{p+1}{2}+\frac{p+3}{2}.
\end{eqnarray*}
Here the approximate values were obtained by using $N_q^{(p)}+1 \approx p/2$. We can see clearly the estimates for $M, F$ scale as $\frac{p!}{2^p}$ which is unfavorable for high order methods.

This motivates the new scheme in this paper by
replacing Gaussian quadrature with Gregory quadrature, which retains high-order accuracy but enables the use of equispaced nodes for computational efficiency. This, in turn, allows for reuse of computations resulting in linear complexity for $M$ and $F$. The details are described in the next section.

\section{Gregory Nested Picard Iteration Schemes}\label{sec:Gregory_NPI_Schemes}
We now describe the Gregory quadrature based NPI scheme for the computation of a $p$th-order approximation. 
We consider discretization of the Duhamel formula
\begin{multline*} \label{eq:Duhamel_form_GREG}
\rho(t^n+\Delta T)
=
U(t^n+\Delta T,t^n)\rho(t^n)U^\dagger(t^n+\Delta T,t^n)
\\ +
\int_{t^n}^{t^n+\Delta T}
U(t^n+\Delta T,s)\,\mathcal{L}_L\rho(s)\,U^\dagger(t^n+\Delta T,s)\,ds.
\end{multline*}
Assuming that the density matrix has been approximated at desired order of accuracy at times $t^n,t^{n+1},\ldots,t^{n+N_q^{(p)}-1},$ we then compute the density matrix at time $t^{n+N_q^{(p)}}$ by the following NPI 
\begin{equation}\label{eq:NPI_general_new_alt}
\begin{aligned}
&\rho^{n+\Nnode^{(p)},k+1} =  U^{(p)}(t^{n+\Nnode^{(p)}},t^n)\rho^n U^{(p)}(t^{n+\Nnode^{(p)}},t^n)^\dagger \\
& + \Delta t\sum_{j=0}^{\Nnode^{(p)}}\omega_j^{(p)} U^{(p)}(t^{n+\Nnode^{(p)}},t^{n+j})\mathcal{L}_L\rho^{n+j,k} U^{(p)} (t^{n+\Nnode^{(p)}},t^{n+j})^\dagger,
\end{aligned}
\end{equation}
for $k=0,1,\dots,p-1$. 
The iteration is started by using the approximation at the previous time step, that is,
\[
\rho^{n+\Nnode^{(p)},0}=\rho^{n+\Nnode^{(p)}-1}.
\]
Note that in the iteration $\rho^{n+j,k} := \rho^{n+j}, j = 0,\ldots, \Nnode^{(p)}-1$ for all $k$.

This scheme is different from the NPI scheme from \cite{hu2025arbitrary} in multiple ways. First, the weights,  $\omega_j^{(p)},$ are now the Gregory quadrature weights of order $p$ (which are positive up to order 9, see Table~\ref{table:gregory}). Second, here the flow operators are always applied at the order of the scheme, while they were applied at successively increasing order in \cite{hu2025arbitrary}. Third, and from a computational efficiency point of view this is perhaps the most important improvement, the formulation allows us to split the quadrature term into  
\begin{multline*}
\Delta t\omega_{\Nnode^{(p)}}^{(p)} \boxed{ \mathcal{L}_L\rho^{n+\Nnode^{(p)},k}}\\
+ \Delta t\sum_{j=0}^{\Nnode^{(p)}-1}\omega_j^{(p)}
U^{(p)}(t^{n+\Nnode^{(p)}},t^{n+j})\mathcal{L}_L\rho^{n+j, k}
U^{(p)}(t^{n+\Nnode^{(p)}},t^{n+j})^\dagger.
\end{multline*}
Here the key observation is that we let $\rho^{n+j, k}=\rho^{n+j},$ therefore throughout the iteration on $k,$ the boxed expression is the \emph{only term that needs to be updated}, and it \emph{does not require} the application of a flow operator. Introducing $\Tilde{\rho}^{n+\Nnode^{(p)}}$ as 
\begin{multline*}
\Tilde{\rho}^{n+\Nnode^{(p)}} :=U^{(p)}(t^{n+\Nnode^{(p)}},t^n)\rho^n U^{(p)}(t^{n+\Nnode^{(p)}},t^n)^\dagger \\
+\Delta t\sum_{j=0}^{\Nnode^{(p)}-1}\omega_j^{(p)}
U^{(p)}(t^{n+\Nnode^{(p)}},t^{n+j})\mathcal{L}_L\rho^{n+j}
U^{(p)}(t^{n+\Nnode^{(p)}},t^{n+j})^\dagger
\end{multline*}
the Gregory-NPI scheme can thus be written compactly as
\begin{equation*}
\rho^{n+\Nnode^{(p)},k+1} := \Tilde{\rho }^{n+\Nnode^{(p)}}+ \Delta t\omega_{\Nnode^{(p)}}^{(p)}  \mathcal{L}_L\rho^{n+\Nnode^{(p)},k}, \ \ k = 0,\ldots,p-1.
\end{equation*}

As just mentioned, $\tilde{\rho}^{\,n+\Nnode^{(p)}}$ can be precomputed before the iteration starts, and that no flow operators are involved during the subsequent iteration stages. The algorithm for Gregory-NPI is described in Algorithm~\ref{alg:algorithm_NPI_Gregory_rho_compact}. 

More generally, to obtain a solution of order $p$, the proposed scheme requires only $(N_q^{(p)}+1)$ flow-operator evaluations. This is a substantial improvement over the Gaussian-quadrature-based NPI schemes.

\begin{algorithm}[]
\caption{Nested Picard iteration scheme with Gregory quadrature, as defined in \eqref{eq:NPI_general_new_alt}}
\label{alg:algorithm_NPI_Gregory_rho_compact}
\begin{algorithmic}[1]
\Require $\rho^n,\ldots,\rho^{n+\Nnode^{(p)}-1}$, $\Delta t$, order $p$, weights $\{\omega_j^{(p)}\}_{j=0}^{\Nnode^{(p)}}$
\Ensure $\rho^{n+\Nnode^{(p)}}$

\State Form
\[
  \tilde{\rho}^{\,n+\Nnode^{(p)}} =
    U^{(p)}(t^{n+\Nnode^{(p)}},t^n)\rho^n
    U^{(p)}(t^{n+\Nnode^{(p)}},t^n)^\dagger
\]
\[
  \qquad
 + \Delta t \sum_{j=0}^{\Nnode^{(p)}-1}\omega_j^{(p)}
 U^{(p)}(t^{n+\Nnode^{(p)}},t^{n+j})
 \mathcal{L}_L\rho^{n+j}
 U^{(p)}(t^{n+\Nnode^{(p)}},t^{n+j})^\dagger .
\]

    \State Set $\rho^{n+\Nnode^{(p)},0} \gets \rho^{\,n+\Nnode^{(p)}-1}$

    \For{$k=1,\dots,p$}
        \State
        \[
        \rho^{n+\Nnode^{(p)},k}
        =
        \tilde{\rho}^{\,n+\Nnode^{(p)}}
        + \Delta t\,\omega_{\Nnode^{(p)}}^{(p)}
        \mathcal{L}_L\rho^{n+\Nnode^{(p)},k-1}
        \]
    \EndFor

\State \Return $\rho^{n+\Nnode^{(p)},p}$
\end{algorithmic}
\end{algorithm}

We note that the Gregory NPI method belongs to the general linear methods (GLM) \cite{jackiewicz2009general}. The first few steps
\[
\rho^1,\dots,\rho^{\Nnode^{(p)}-1}
\]
need to be initialized. In practice, one often uses either a higher-order one-step method or a finer time step size in this startup phase to ensure that the starting error does not contaminate the overall convergence order. Here we generate these starting values on an auxiliary fine grid with step size $\Delta t/m$, where $m=O(\Delta t^{-(p-3)/3})$ so that $(\Delta t/m)^3=O(\Delta t^p)$. The second-order scheme is first used to compute the approximations at the first two fine-grid nodes, not counting the initial value. The order is then increased successively: for $q=3,\ldots,p-1$, the order-$q$ scheme is used to compute the next two values, and the order-$p$ scheme is applied from the fine-grid $(2p-3)$-th node onward. Other approaches are of course also possible.

The approximation of the flow operator $U(t^{n+k},t^n)$ associated with
\[
\frac{dV(t)}{dt}=J(t)V(t), \qquad V(t^n)=V_n.
\]
have been discussed in prior work \cite{hu2025arbitrary}.
In our numerical simulations, explicit approximations of the flow operator are constructed using Runge–Kutta methods. For implicit approximations, we use the trapezoidal rule for the second-order scheme and the method proposed in \cite{puzynin2000magnus} for the fourth-order scheme. Other high-order implicit time schemes can be used but we do not consider those in this paper. 

\section{Stability Analysis}\label{sec: stability}
Theoretical results of GLM \cite{jackiewicz2009general} have been discussed in the literature. The proposed Gregory NPI methods are not designed to be a GLM for general ODEs, but they are tailored for the Lindblad equation. Therefore,  in this paper, we focus the analysis   specifically on the Lindblad equation. 

Although trace renormalization is a nonlinear operation that keeps the density matrix bounded, it may compromise the accuracy of the scheme. This effect is illustrated in the Qudit-cavity example of Section~\ref{sec: experiments: Qudit cavity}, where the numerical solution remains bounded but is not  accurate. It is therefore necessary to study the stability of \eqref{eq:NPI_general_new_alt} without trace renormalization, which is what we will do in this section.

Similar to \cite{hu2025arbitrary}, we consider a single-qubit test equation with 
the Hamiltonian and the Lindblad operators of the model  are given by (note that we use $\Omega$ instead of $\omega$ so as to not confuse it with the quadrature weights)
\[
H =  \Omega \bm{0&0\\0&1}, \mathcal{L}_1= \frac{1}{\sqrt{T_1}}\bm{0&1\\0&0},  \ \ \mathcal{L}_2 = \frac{1}{\sqrt{T_2}}\bm{0&0\\0&1},
\]
This will give rise to
\[
J = -iH -\frac{1}{2} \mathcal{L}_1^\dagger \mathcal{L}_1 -\frac{1}{2} \mathcal{L}_2^\dagger \mathcal{L}_2 = \bm{0&0\\0&-i\Omega  -\frac{1}{2T_1}-\frac{1}{2T_2}}.
\]

 Denote   $\mathcal{W} := \Delta t \Omega$, $\mathcal{T}_1 := \frac{\Delta t}{T_1} $ and $\mathcal{T}_2 := \frac{\Delta t}{T_2} $. All approximations of the flow operator $U(t^{n+k},t^n)$ considered in the following subsections take the form

\[
\bm{1&0\\0&\varsigma^{(p,j)}(\mathcal{W},\mathcal{T}_1, \mathcal{T}_2, k)},
\]
where $k$ denotes the number of time steps from $t^n$ to $t^{n+k}$, $p$ denotes the order of approximation, and $j\in\{0,1\}$ indicates whether the approximation is explicit or implicit. More precisely,
\[
j=
\begin{cases}
0, & \text{if the approximation is explicit},\\
1, & \text{if the approximation is implicit}.
\end{cases}
\]

The quantity $\varsigma^{(p,j)}(\mathcal{W},\mathcal{T}_1,\mathcal{T}_2,k)$ depends on the particular approximation of the flow. When there is no ambiguity,  $\varsigma^{(p,j)}(k)$ means $\varsigma^{(p,j)}(\mathcal{W},\mathcal{T}_1, \mathcal{T}_2, k)$. We also use $\varsigma^{(p)}(k)$ when whether the flow operator is implicit or explicit is not specified.

We also denote $A(k):= -\frac{k}{2}\frac{\Delta t}{T_1} - \frac{k}{2}\frac{\Delta t}{T_2}-ik\Omega\Delta t  = -\frac{k}{2}\mathcal{T}_1-\frac{k}{2}\mathcal{T}_2-ik\mathcal{W} $, which will be used in the stability analysis. Next, we derive matrix representations for the linear operators used in our numerical scheme, which will be needed in the stability analysis. Using columnwise vectorization and the identity
\[
\operatorname{vec}(AXB)=(B^T\otimes A)\operatorname{vec}(X),
\]
the map \(\rho\mapsto U\rho U^\dagger\) is represented by
\[
\operatorname{vec}(U\rho U^\dagger)
=
(\overline{U}\otimes U)\operatorname{vec}(\rho),
\]
where \(\overline{U}\) denotes the entrywise complex conjugate of \(U\). It follows that the matrices
\begin{align*}
P_k^{(p)}&=\bm{1&0&0&0\\0&\overline{\varsigma^{(p)}(k)}&0&0\\0&0&\varsigma^{(p)}(k)&0\\0&0&0&|\varsigma^{(p)}(k)|^2},\\[1mm]
Q_{U^{(p)}L,k} &=
\bm{0&0&0&\frac{1}{T_1}\\0&0&0&0\\0&0&0&0\\0&0&0&\frac{|\varsigma^{(p)}(k)|^2}{T_2}},\\[1mm]
Q_L &=
\bm{0&0&0&\frac{1}{T_1}\\0&0&0&0\\0&0&0&0\\0&0&0&\frac{1}{T_2}},
\end{align*}
represent the operators
\begin{align*}
   & \operatorname{vec}(\rho) \mapsto \operatorname{vec}\left(U^{(p)}(t^{n+k},t^n)\rho \bigl(U^{(p)}(t^{n+k},t^n)\bigr)^\dagger\right),\\
&\operatorname{vec}(\rho) \mapsto \operatorname{vec}\left(U^{(p)}(t^{n+k},t^n)\mathcal{L}_L\rho \bigl(U^{(p)}(t^{n+k},t^n)\bigr)^\dagger\right),\\
&\operatorname{vec}(\rho) \mapsto \operatorname{vec}\left(\mathcal{L}_L\rho\right),
\end{align*}
respectively. These three operators form the main components of our schemes. Their matrix representations indicate that $\varsigma^{(p)}(k)$, $T_1$, and $T_2$ play essential roles in determining stability.

We begin with the stability analysis of the second-order scheme, which is essentially a one-step method. Similar analysis has been carried out in \cite{hu2025arbitrary}. We then proceed to the higher-order schemes, whose stability must be analyzed in the framework of GLM.
\subsection{Second-Order Scheme}
The second-order approximation of the solution is defined by
\begin{align}
\rho^{n+1,0} & =  \rho^n \nonumber\\
\rho^{n+1,1} &=  U^{(2)}(t^{n+1},t^n) \rho^n U^{(2)}(t^{n+1},t^n) ^\dagger +\nonumber \\ 
&\qquad {\Delta t}\Big(\omega^{(2)}_0  U^{(2)}(t^{n+1},t^{n}) \mathcal{L}_L\rho^n U^{(2)}(t^{n+1},t^{n}) ^\dagger 
 +\omega^{(2)}_1 \mathcal{L}_L {\rho^{n+1,0}} \Big),\nonumber\\
\rho^{n+1,2} &=  U^{(2)}(t^{n+1},t^n) \rho^n U^{(2)}(t^{n+1},t^n) ^\dagger +\nonumber\\
&\qquad {\Delta t}\Big(\omega^{(2)}_0  U^{(2)}(t^{n+1},t^{n}) \mathcal{L}_L\rho^n U^{(2)}(t^{n+1},t^{n}) ^\dagger +\omega^{(2)}_1 \mathcal{L}_L {\rho^{n+1,1}} \Big),\label{eqn: NPI_second-order}
\end{align}
where  $\omega^{(2)}_0 = \omega^{(2)}_1 = \frac{1}{2}$.
For the second-order approximation of the flow operator, we apply the explicit midpoint method in the explicit case and the trapezoidal rule in the implicit case, i.e.
\[
U^{(2)}(t^{n+k},t^n)=
\begin{cases}
I+k\Delta t\,J+\dfrac{(k\Delta t)^2}{2}J^2, & \text{if } j=0, \qquad{\rm i.e., explicit}, \\[1.2ex]
\left(I-\dfrac{k\Delta t}{2}J\right)^{-1}\left(I+\dfrac{k\Delta t}{2}J\right), & \text{if } j=1, \qquad{\rm i.e., implicit}.
\end{cases}
\]
With this flow operator we get
\[
\varsigma^{(2,j)}(k)=
\begin{cases}
1+A(k)
+\dfrac{1}{2}A(k)^2, & \text{if } j=0, \\[2ex]
\dfrac{
1+\dfrac{1}{2}A(k)
}{
1-\dfrac{1}{2}A(k)
}, & \text{if } j=1.
\end{cases}
\]

To analyze the stability of the scheme, we first rewrite the matrix equation in vectorized form. This yields an amplification matrix representation of the method, allowing the stability analysis to be reduced to a spectral study of the corresponding matrix, in particular through its eigenvalues and spectral radius.
Throughout the following discussion, boldface symbols denote vectorized matrices. In particular, $\boldsymbol{\rho}$ denotes the vectorized form of $\rho$. Vectorizing \eqref{eqn: NPI_second-order} yields
\begin{align*}
   \boldsymbol{ \rho}^{n+1,1}  & =P_1^{(2)} \boldsymbol{\rho}^n  + \Delta t\omega_0^{(2)} Q_{U^{(2)}L,1} \boldsymbol{\rho}^n  + \Delta t\omega_1^{(2)} Q_L {\boldsymbol{\rho}^{n+1,0}} , \\
    \boldsymbol{\rho}^{n+1,2} & = P_1^{(2)} \boldsymbol{\rho}^n + \Delta t\omega_0^{(2)} Q_{U^{(2)}L,1} \boldsymbol{\rho}^n + \Delta t\omega_1^{(2)} Q_L {\boldsymbol{\rho}^{n+1,1}}.
\end{align*}

\begin{figure}[htb]
\centering
\includegraphics[width=0.3\textwidth, trim=75mm 64mm 76mm 66mm, clip=true]{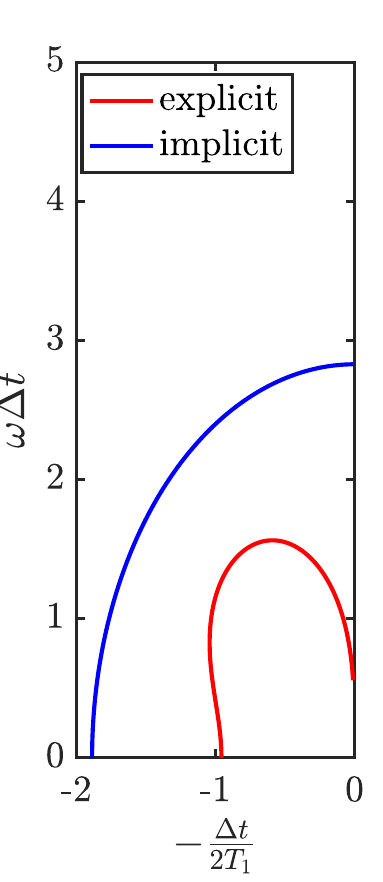}
\includegraphics[width=0.3\textwidth, trim=75mm 64mm 76mm 66mm, clip=true]{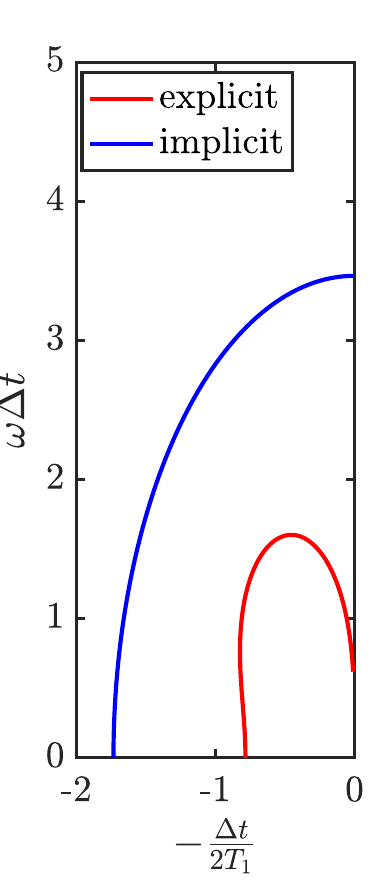}
\includegraphics[width=0.3\textwidth, trim=75mm 64mm 76mm 66mm, clip=true]{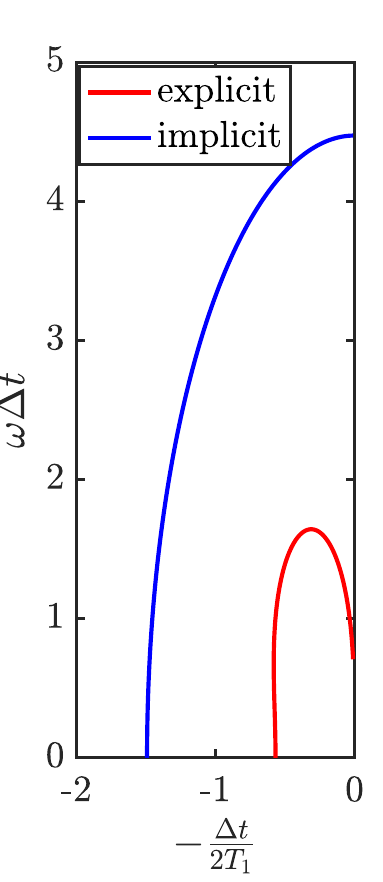}
\caption{ Stability regions of the second-order schemes with explicit and implicit flow operators, for  $T_2 = 2T_1$ (left), $T_2 = T_1$ (middle) and $T_2 = T_1/2$ (right).\label{fig:stability_2}}
\end{figure}

Substituting the expression for $\boldsymbol{\rho}^{n+1,1}$ into the equation for $\boldsymbol{\rho}^{n+1,2}$ yields
\begin{multline*}
\boldsymbol{\rho}^{n+1,2} = \Bigl( P_1^{(2)} + \Delta t\omega_0^{(2)} Q_{U^{(2)}L,1} + {\Delta t}\omega_1^{(2)} Q_L  P_1^{(2)} +\\(\Delta t)^2 \omega_0^{(2)} \omega_1^{(2)} Q_L Q_{U^{(2)}L,1}+ (\Delta t)^2(\omega_1^{(2)} )^2( Q_L )^2
 \Bigr)\boldsymbol{\rho}^n.     
\end{multline*}
The amplification matrix of the second-order approximation is 
\begin{align*}
& P_1^{(2)} + \Delta t\omega_0^{(2)} Q_{U^{(2)}L,1} + {\Delta t}\omega_1^{(2)} Q_L  P_1^{(2)} +(\Delta t)^2 \omega_0^{(2)} \omega_1^{(2)} Q_L Q_{U^{(2)}L,1}\\
&+ (\Delta t)^2(\omega_1^{(2)} )^2( Q_L )^2 =\\
& \bm{1&0&0&
| \varsigma^{(2)}(1)|^2 \left(\frac{\Delta t}{T_1}\omega_1^{(2)} +\omega_0^{(2)}\omega_1^{(2)} \frac{(\Delta t)^2}{T_1T_2}  \right)
+\frac{\Delta t}{T_1}(\omega_0^{(2)}
+  \frac{\omega_1^{(2)} \Delta t}{T_2})
\\
0&\overline{\varsigma^{(2)}(1)}&0&0
\\
0&0&\varsigma^{(2)}(1)&0
\\
0&0&0&|\varsigma^{(2)}(1)|^2 
\left( 1+ \frac{\Delta t}{T_2} 
+  \omega_0^{(2)}\omega_1^{(2)} (\frac{\Delta t}{T_2})^2 \right)
+  (\frac{\omega_1^{(2)}\Delta t}{T_2})^2}.
\end{align*}
So the stability condition for the second-order scheme is 
\[
\left| |\varsigma^{(2)}(1)|^2 \left( 1+ \frac{\Delta t}{T_2} + \omega_0^{(2)}\omega_1^{(2)} (\frac{\Delta t}{T_2})^2 \right)+ (\frac{\omega_1^{(2)}\Delta t}{T_2})^2 \right|  \le 1.
\]
Note that, as $\Delta t/T_2 \to 0$, the stability condition approaches
\[
\left| \varsigma^{(2)}(1)\right|\le 1,
\]
that is, the stability condition associated with the second-order flow solver. For illustration, Figure~\ref{fig:stability_2} shows the stability regions of the second-order schemes for three representative relations between the pure dephasing time and the decay time, namely,
\[
T_2=2T_1,\qquad T_2=T_1,\qquad T_2=\frac{T_1}{2}.
\]
In all three cases, the scheme with the implicit flow solver has a much larger stability region than the corresponding scheme with the explicit flow solver.

\subsection{Fourth-Order Scheme}
The fourth-order scheme is based on the fourth-order Gregory quadrature rule, which requires at least six quadrature nodes. The fourth-order approximation is defined by
\begin{align}
\rho^{n+5,0} & =  \rho^{n+4} \nonumber\\
\rho^{n+5,k} &=  U^{(4)}(t^{n+5},t^n) \rho^n U^{(4)}(t^{n+5},t^n) ^\dagger +\\
&\Delta t\Big( \omega^{(4)}_0 U^{(4)}(t^{n+5},t^n) \mathcal{L}_L \rho^n U^{(4)}(t^{n+5},t^n) ^\dagger+\nonumber\\
& \omega^{(4)}_1 U^{(4)}(t^{n+5},t^{n+1}) \mathcal{L}_L \rho^{n+1} U^{(4)}(t^{n+5},t^{n+1}) ^\dagger+\nonumber\\
&
\omega^{(4)}_2 U^{(4)}(t^{n+5},t^{n+2}) \mathcal{L}_L \rho^{n+2} U^{(4)}(t^{n+5},t^{n+2}) ^\dagger+\nonumber\\
&\omega^{(4)}_3 U^{(4)}(t^{n+5},t^{n+3}) \mathcal{L}_L \rho^{n+3} U^{(4)}(t^{n+5},t^{n+3}) ^\dagger+\nonumber\\
&
\omega^{(4)}_4U^{(4)}(t^{n+5},t^{n+4}) \mathcal{L}_L \rho^{n+4} U^{(4)}(t^{n+5},t^{n+4}) ^\dagger \Big)\nonumber\\
& + \Delta t \omega^{(4)}_5\mathcal{L}_L { \rho^{n+5,k-1}}, \qquad   \label{eqn: approximation_order_4}
\end{align}
for $k=1,2,3,4.$ Here, the Gregory weights are $\omega^{(4)}_0 =\omega^{(4)}_5 = \frac{3}{8}, \omega^{(4)}_1 =\omega^{(4)}_4 = \frac{7}{6}, \omega^{(4)}_2 =\omega^{(4)}_3 = \frac{23}{24}$. 
For the fourth-order approximation of the flow operator, we can apply the classical fourth-order  Runge-Kutta  method in the explicit case, namely, 
\[
U^{(4)}(t^{n+k},t^n)=
I + k\Delta t J + \frac{(k\Delta t)^2}{2}J^2+\frac{(k\Delta t)^3}{6}J^3+\frac{(k\Delta t)^4}{24}J^4,
\]
or the fourth-order  implicit scheme from \cite{puzynin2000magnus} in the implicit case, namely,
\[
U^{(4)}(t^{n+k},t^n)=
(I +  i\frac{k\Delta t}{4}\overline{d} J)^{-1}(I + i\frac{k\Delta t}{4} d  J)\,
(I - i\frac{k\Delta t}{4}d J)^{-1}(I - i\frac{k\Delta t}{4}\overline{d}  J),
\] 
with $d=\frac{1}{\sqrt{3}}-i$. It follows that
\[
\varsigma^{(4,j)}(k)=\begin{cases}
    1+ A(k)+\frac{1}{2}A(k)^2 +\frac{1}{6}A(k)^3
  +\frac{1}{24}A(k)^4,& \text{if } j=0,\\
   \dfrac{( 1+i\frac{d}{4}A(k))(1-i\frac{\bar{d}}{4}A(k)) }{ (1+i\frac{\bar{d}}{4}A(k)) ( 1-i\frac{d}{4}A(k))},  & \text{if } j=1.
\end{cases}
\]
 We present the following stability result whose proof can be found in the Appendix.
\begin{theorem}\label{thm: order 4}
   The fourth-order scheme is stable if all roots of the fifth-degree polynomial
\[
\lambda^5
-\sum_{k=1}^{5}\lambda^{5-k}\omega^{(4)}_{5-k}\,B^{(4)}\frac{\Delta t}{T_2}\left|\varsigma^{(4)}(k)\right|^2
-\lambda^4(\omega^{(4)}_5)^4\left(\frac{\Delta t}{T_2}\right)^4
-B^{(4)}\left|\varsigma^{(4)}(5)\right|^2
\]
lie in the unit disk, and
\[
\left|\varsigma^{(4)}(5)\right|\le1,
\]
where  $B^{(4)}= \sum_{k=0}^3 (\omega^{(4)}_5)^k(\frac{\Delta t}{T_2})^k
$. 
\end{theorem}
Note that, for sufficiently large dephasing time, or equivalently, for sufficiently small $\frac{\Delta t}{T_2}$, the stability region approaches the one  defined by
\[
\left|\varsigma^{(4)}(5)\right|\le1.
\]
\begin{figure}[]
\centering
\includegraphics[width=0.35\textwidth, trim=46mm 64mm 46mm 66mm, clip=true]{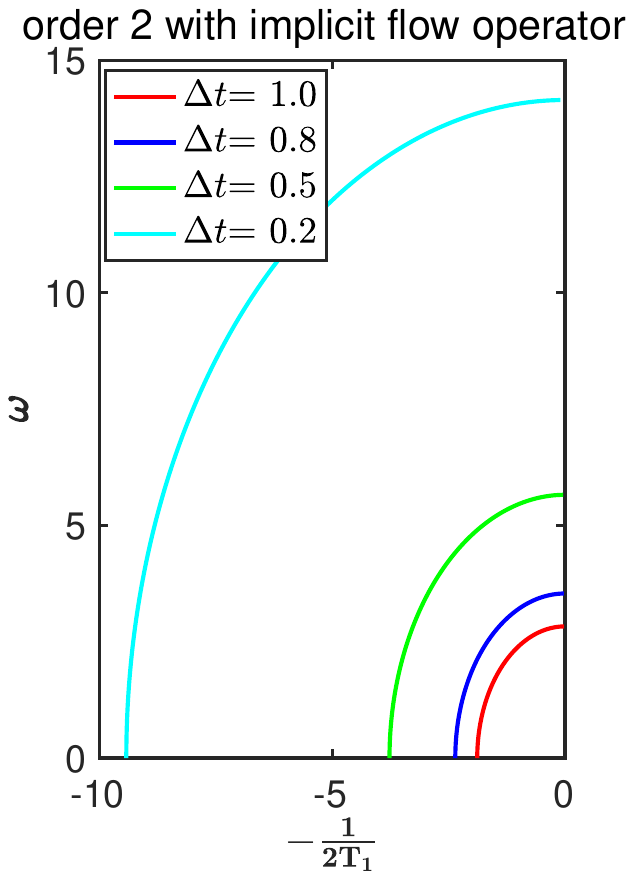}
\includegraphics[width=0.35\textwidth, trim=46mm 64mm 46mm 66mm, clip=true]{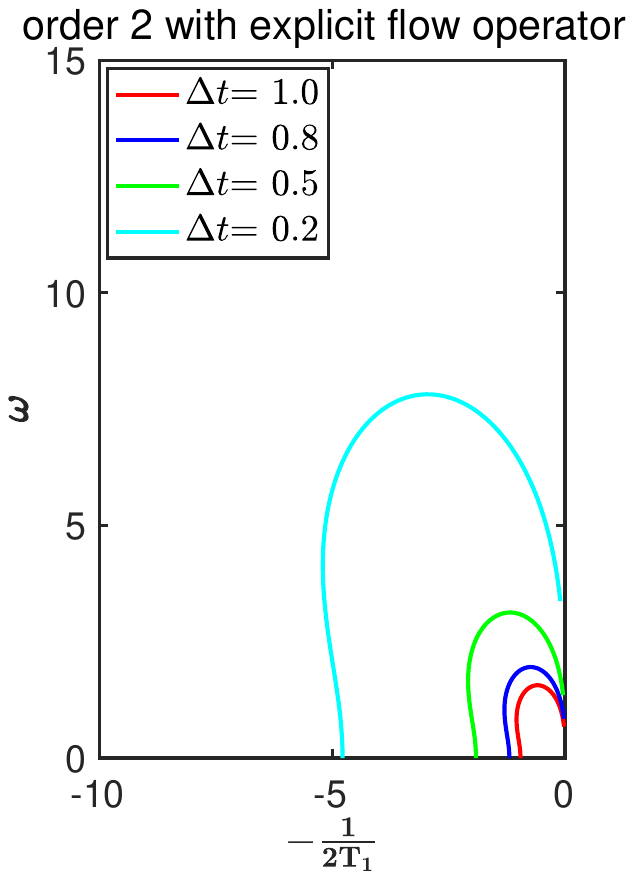}
\includegraphics[width=0.35\textwidth, trim=46mm 64mm 46mm 66mm, clip=true]{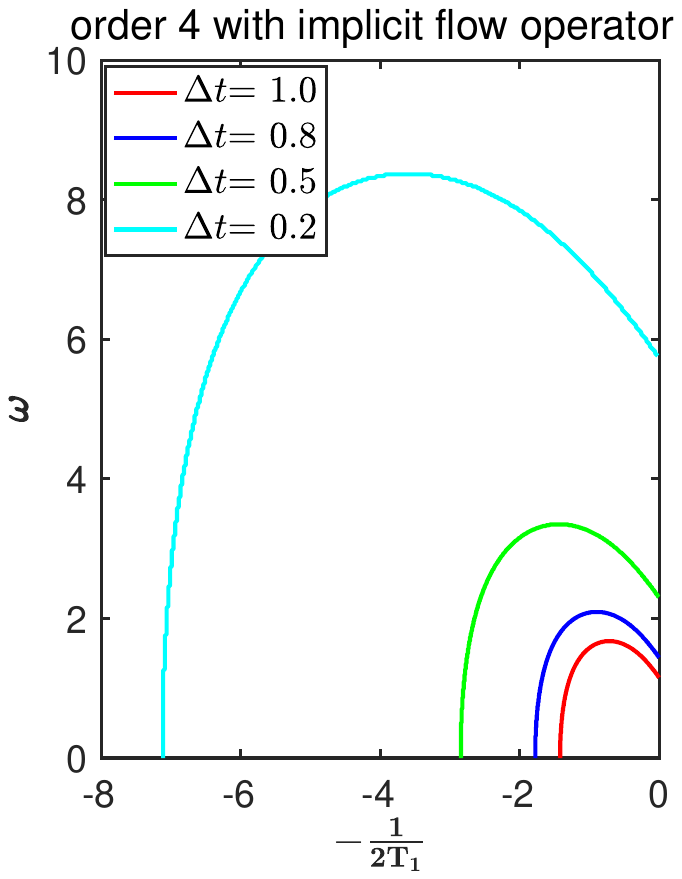}
\includegraphics[width=0.35\textwidth, trim=46mm 64mm 46mm 66mm, clip=true]{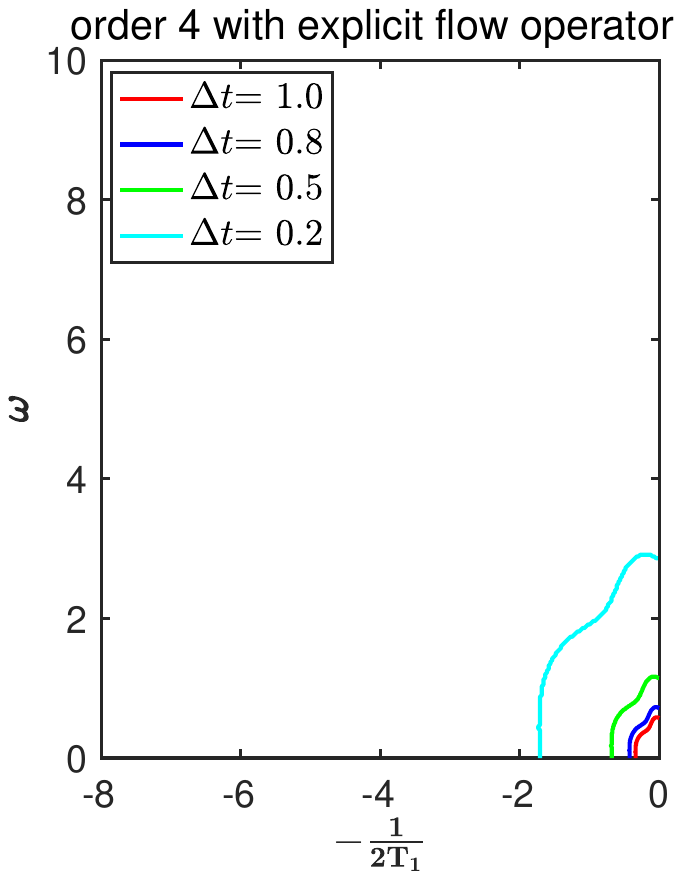}
\caption{ Stability regions of the second-order and fourth-order schemes with explicit and implicit flow operators, for  $T_2 = 2T_1$, fine time-step sizes $\Delta t = 0.2$, $0.5$, $0.8$, and $1.0$. }
\label{fig:stability_2_4}
\end{figure}

In other words, when the dephasing time is large, the stability region is determined primarily by the stability of the flow operator. This phenomenon is also observed for the second-order method.
In Figure~\ref{fig:stability_2_4}, we plot the stability regions of the second- and fourth-order schemes for time step sizes $\Delta t=0.2, 0.5, 0.8, 1.0$, under the constraint $T_2=2T_1$. For current quantum devices, the relaxation and dephasing times are typically large enough that $1/(2T_1)$ remains small, and practical choices of $\Delta t$ therefore fall within the stability region.

\subsection{General $p$-th Order Scheme}
We can generalize the stability results for higher order. 
\begin{theorem}\label{thm: order p}
   The $p$-th order scheme is stable if all roots of the $(2p-3)$-th-degree polynomial
\begin{align*}
    \lambda^{2p-3}
-\sum_{k=1}^{2p-3}\lambda^{2p-3-k}\omega^{(p)}_{2p-3-k}\,B^{(p)}\frac{\Delta t}{T_2}\left|\varsigma^{(p)}(k)\right|^2
-&\\
\lambda^{2p-4}(\omega^{(p)}_{2p-3})^{p}\left(\frac{\Delta t}{T_2}\right)^p
-B^{(p)}\left|\varsigma^{(p)}(2p-3)\right|^2&
\end{align*}
lie in the  unit disk, and
\[
\left|\varsigma^{(p)}(2p-3)\right|\le1,
\]
where  $B^{(p)}= \sum_{k=0}^{p-1} (\omega^{(p)}_{2p-3})^k(\frac{\Delta t}{T_2})^k
$. 
\end{theorem}
The proof is similar to the proof of Theorem \ref{thm: order 4} and is skipped for brevity.

Here, similar to previous cases, when $\frac{\Delta t}{T_2}\rightarrow 0$, i.e., when the pure dephasing time is sufficiently large, the stability condition is governed almost entirely by the approximation of the flow operator. Consequently, for higher-order schemes in which the flow operator is approximated by a Runge–Kutta method, the stability of the overall scheme is essentially inherited from that of the underlying Runge–Kutta method.

\section{Low-Rank Truncation}\label{sec:truncation}
We now describe the low-rank truncation procedure for scheme  \eqref{eq:NPI_general_new_alt}. We represent the density matrix in factored form as
\[
\rho(t)=V(t)V(t)^\dagger,
\]
where \(\rho(t)\in\mathbb{C}^{d\times d}\) and \(V(t)\in\mathbb{C}^{d\times r(t)}\) with \(r(t)\ll d\). Thus, instead of storing and evolving \(\rho(t)\) directly, we work with its low-rank factor \(V(t)\).

At each update, the truncation is performed in two stages. The first stage is applied after the offline step. We write the intermediate approximation \(\widetilde{\rho}^{\,n+\Nnode^{(p)}}\) as
\[
\widetilde{\rho}^{\,n+\Nnode^{(p)}}
=
\widetilde{V}^{\,n+\Nnode^{(p)}}
\bigl(\widetilde{V}^{\,n+\Nnode^{(p)}}\bigr)^\dagger,
\]
where
\[
\widetilde{V}^{\,n+\Nnode^{(p)}}
=
\Bigl[
U^{(p)}(t^{n+\Nnode^{(p)}},t^n)V^n
~~ \cdots ~~
\sqrt{\omega_j^{(p)}\Delta t}\,
U^{(p)}(t^{n+\Nnode^{(p)}},t^{n+j})\mathcal{L}_{\alpha}V^{n+j}
\quad \cdots
\Bigr].
\]
Here, $ j=0,\dots,\Nnode^{(p)}-1.$
The second stage is applied during the nested Picard iteration. We set
\[
\rho^{n+\Nnode^{(p)},0}
=
V^{n+\Nnode^{(p)}-1}\bigl(V^{n+\Nnode^{(p)}-1}\bigr)^\dagger,
\]
and, for \(k=1,\dots,p\),
\[
\rho^{n+\Nnode^{(p)},k}
=
V^{n+\Nnode^{(p)},k}\bigl(V^{n+\Nnode^{(p)},k}\bigr)^\dagger,
\]
where
\[
V^{n+\Nnode^{(p)},k}
=
\Bigl[
\widetilde{V}^{\,n+\Nnode^{(p)}}
\qquad
\sqrt{\omega_{\Nnode}^{(p)}\Delta t}\,
\mathcal{L}_{\alpha}V^{n+\Nnode^{(p)},k-1}
\Bigr].
\]
To maintain a low-rank representation, we truncate
\[
\widetilde{V}^{\,n+\Nnode^{(p)}},\qquad
V^{n+\Nnode^{(p)},k},\quad k=1,\dots,p,
\]
by means of a truncated singular value decomposition. Following \cite{hu2025arbitrary}, we choose the truncation tolerance proportional to the local truncation error. For a method of order \(p\), we take
\[
\varepsilon_{\mathrm{tr}}=(\kappa\Delta t)^{p+1},
\]
where \(\kappa>0\) is a prescribed constant.

When the number of columns in \(V^{n+\Nnode,k}\) becomes large, a direct truncated SVD may be computationally expensive. In that case, we use the same hierarchical truncation procedure as in \cite{hu2025arbitrary}, namely, we group the columns and apply truncated SVD successively to each group. Since each truncation step is performed at the level of the factors \(\widetilde{V}^{\,n+\Nnode}\) and \(V^{n+\Nnode,k}\), the resulting density matrix after truncation remains completely positive.

After advancing to time \(t^{n+\Nnode^{(p)}}\), we apply the trace renormalization
\[
\rho^{n+\Nnode^{(p)}}
\leftarrow
\frac{\rho^{n+\Nnode^{(p)}}}{\operatorname{Tr}(\rho^{n+\Nnode^{(p)}})}.
\]
We note that the dominant cost of the truncation procedure arises in the offline stage, while the online stage is relatively inexpensive. Compared with the scheme in \cite{hu2025arbitrary}, the truncation procedure in the present work is computationally cheaper.

\section{Numerical Experiments}\label{sec: numerical experiments}
Direct simulation in the laboratory frame becomes expensive when $H(t)$ contains highly oscillatory terms, since these induce stiffness and force very small time steps. To reduce this difficulty, we work in a rotating frame and apply the rotating-wave approximation, thereby eliminating much of the rapidly oscillatory behavior.

Given rotation frequencies $\omega_k^r$, the transformed drift and control Hamiltonians are
\begin{align*}
\widetilde{H}_d(t)
=&
\sum_{k=0}^{Q-1}
\left(
(\omega_k-\omega_k^r)a_k^\dagger a_k
-\frac{\xi_k}{2}\, a_k^\dagger a_k^\dagger a_k a_k
-\sum_{l>k}\xi_{kl} a_k^\dagger a_k a_l^\dagger a_l
\right) \\
&\quad
+\sum_{k=0}^{Q-1}\sum_{l>k}
J_{kl}
\left(
\cos(\eta_{kl}t)(a_k^\dagger a_l+a_k a_l^\dagger)
+i\sin(\eta_{kl}t)(a_k^\dagger a_l-a_k a_l^\dagger)
\right), \\
\widetilde{H}_c(t)
=&
\sum_{k=0}^{Q-1}
\left(
p^k(\vec{\alpha}^k,t)(a_k+a_k^\dagger)
+i q^k(\vec{\alpha}^k,t)(a_k-a_k^\dagger)
\right),
\end{align*}
where
\[
\eta_{kl}:=\omega_k^r-\omega_l^r.
\]
The laboratory-frame control can then be approximated by
\[
f^k(t)\approx 2p^k(t)\cos(\omega_k^r t)-2q^k(t)\sin(\omega_k^r t).
\]

The relaxation and dephasing parameters used in the following experiments are chosen to be consistent with experimentally reported coherence scales in superconducting circuits. Recent transmon devices report relaxation and coherence times exceeding $100\,\mu\mathrm{s}$, with high-coherence platforms reaching several hundred microseconds, while cavity or resonator memories can exhibit lifetimes from the millisecond regime to several tens of milliseconds. Thus, the values considered here are intended to be representative of experimentally relevant regimes rather than tied to a single device realization.

In the following subsections, we use standard bra--ket notation, where
\[
\ket{j_1j_2\cdots j_Q}:=\ket{j_1}\otimes\ket{j_2}\otimes\cdots\otimes\ket{j_Q}
\]
denotes a tensor-product basis state, and \(\bra{\psi}=\ket{\psi}^\dagger\).

\subsection{Two-Qubit System}
In the first experiment, we investigate the accuracy of several numerical methods for a composite quantum system consisting of two identical qubits with the same transition frequency, that is, $\omega_0=\omega_1$. In the rotating frame, we choose the rotation frequencies to coincide with the physical frequencies, namely, $\omega_k^r=\omega_k$ for $k=0,1$. With this choice, the terms $\sum_k \omega_k a_k^\dagger a_k$ are eliminated from the system Hamiltonian. The dipole-dipole interaction is characterized by a single cross-resonance coefficient $J_{01}=0.2$ (in nondimensionalized units), so that the Hamiltonian reduces to
\[
H = J_{01}(a_0^\dagger a_1 + a_0 a_1^\dagger).
\]
We consider Lindblad operators associated with decay and dephasing processes, with uniform rates
\[
\mu = \frac{1}{T_1^0}=\frac{1}{T_1^1}=\frac{1}{50}, \qquad
\gamma_1=\frac{1}{T_2^0}=\frac{1}{50}, \qquad
\gamma_2=\frac{1}{T_2^1}=\frac{1}{50}.
\]
Define
\[
a=\frac{\gamma_1+\gamma_2}{4}, \qquad
b=\sqrt{4J_{01}^2-a^2}.
\]

For the initial condition $\rho(0)=|10\rangle\langle 10|$, the exact solution at time $t$ is given by
\begin{equation*}
\rho(t) = \begin{pmatrix}
\rho_{11}(t)&0&0&0\\
0&\rho_{22}(t)&\rho_{23}(t)&0\\
0&\rho_{32}(t)&\rho_{33}(t)&0\\
0&0&0&0
\end{pmatrix},
\end{equation*}
where 
\begin{align*}
\rho_{11}(t)&=1-e^{-\mu t},\\    
\rho_{22}(t)
&=
\frac12\left[
e^{-\mu t}
-
e^{-(\mu+a)t}
\left(
\cos(b t)+\frac{a}{b}\sin(b t)
\right)
\right],\\
\rho_{33}(t)
&=
\frac12\left[
e^{-\mu t}
+
e^{-(\mu+a)t}
\left(
\cos(b t)+\frac{a}{b}\sin(b t)
\right)
\right],\\
\rho_{23}(t)
&=
-\,i\,\frac{J_{01}}{b}\,e^{-(\mu+a)t}\sin(b t),
\qquad
\rho_{32}(t)=\,i\,\frac{J_{01}}{b}\,e^{-(\mu+a)t}\sin(b t).
\end{align*}
We simulate the density matrix up to the final time $T=6$ (in nondimensionalized units) using methods of orders two through nine, and measure the error of the density matrix at the final time in the Frobenius norm. Table~\ref{tab:2Q1} reports the errors for the second-, third-, and fourth-order methods. Here, “explicit” and “implicit” indicate whether the flow operator is approximated by an explicit or an implicit method. The observed convergence rates agree with the expected orders, except for the implicit third-order scheme, which exhibits fourth-order convergence. This is likely due to the use of a fourth-order implicit flow operator approximation in the third-order implicit scheme. The same phenomenon was observed in \cite{hu2025arbitrary}.

\begin{figure}[htb]
\centering
\includegraphics[width=0.5\textwidth, trim=9mm 66mm 27mm 74mm, clip=true]{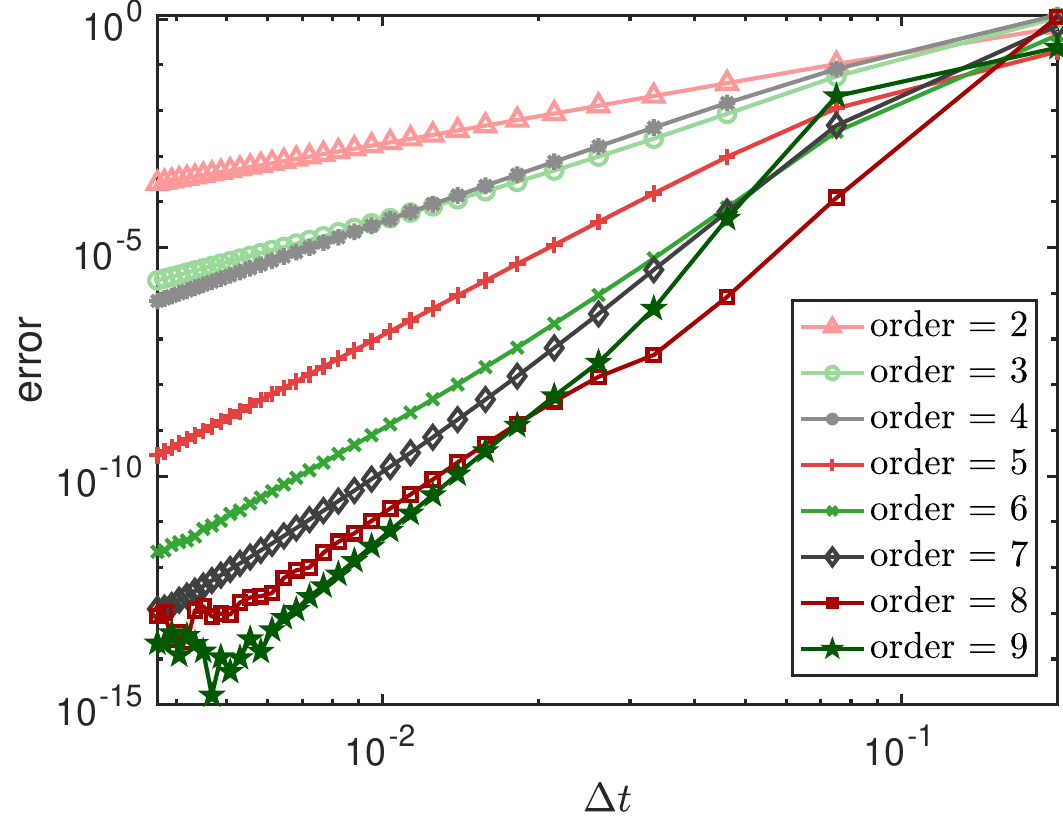}
\caption{Two-Qubit experiment: errors in the Frobenius norm for schemes of orders 2 through 9 with explicit flow operators. 
}
\label{fig:2Q1_err_all_orders}
\end{figure}

Figure~\ref{fig:2Q1_err_all_orders} shows the error versus the time step size for methods of orders two through nine, where the flow operator is approximated by the Runge--Kutta method associated with each scheme. The results again confirm that the observed convergence rates are consistent with the expected orders.

\begin{table}[] 
 \begin{center} 
 \begin{tabular}{| l | c | c | c | c | } 
\hline
\multicolumn{5}{|c|}{second-order explicit} \\  \hline
$N_t$ & 128 & 256 & 512 & 1024 \\\hline
error & 4.10(-2) & 1.03(-2) & 2.57(-3) & 6.44(-4)  \\ \hline 
rate &  & 1.99 & 2.00 & 2.00 \\\hline
\hline
\multicolumn{5}{|c|}{second-order implicit}  \\\hline $N_t$ &128 & 256 & 512 & 1024\\\hline  error & 2.09(-2) & 5.16(-3) & 1.28(-3) & 3.21(-4) \\ \hline rate &  & 2.01 & 2.01 & 2.00 \\ \hline 
\hline
\multicolumn{5}{|c|}{third-order explicit}    \\ \hline
$N_t$ & 96 & 192 & 384 & 768 \\
 \hline 
error & 2.71(-2) & 1.90(-3) & 1.66(-4) & 1.85(-5)  \\ \hline 
rate &  & 3.83 & 3.52 & 3.16 \\\hline 
\hline\multicolumn{5}{|c|}{third-order implicit}\\\hline$N_t$ &96 & 192 & 384 & 768 \\\hline  error & 1.16(-3) & 7.25(-5) & 4.04(-6) & 1.38(-7)\\\hline rate &  & 4.00 & 4.17 & 4.87 \\ \hline 
\hline
\multicolumn{5}{|c|}{fourth-order explicit}  \\ \hline
$N_t$ & 80 & 160 & 320 & 640 \\ \hline  
error & 8.11(-2) & 6.67(-3) & 4.46(-4) & 2.84(-5)  \\ \hline  
rate &  &  3.60 & 3.90 & 3.97 \\\hline\hline
\multicolumn{5}{|c|}{fourth-order implicit} \\ \hline$N_t$ & 80 & 160 & 320 & 640 \\\hline error & 1.73(-2) & 1.17(-3) & 7.48(-5) & 4.73(-6)\\\hline rate &  & 3.89 & 3.96 & 3.99 \\ \hline 
 \end{tabular} 
 \caption{Errors and rates of convergence at the final time for the Two-Qubit  example with a known analytic solution. Here the trace renormalization is used. In this table $4.10(-2) = 4.10 \cdot 10^{-2}$. 
 \label{tab:2Q1}}
 \end{center} 
 \end{table} 

\subsection{Qudit-Resonator System with Control}\label{sec: experiments: Qudit cavity}
In this section, we consider a composite quantum system consisting of a three-level qudit and a twenty-level resonator. The device model is standard and was also employed in \cite{gunther2021quantum} and our earlier work \cite{hu2025arbitrary}.  
We adopt the same parameter values as in \cite{hu2025arbitrary}; these are reported in Table~\ref{table:parameters_QuditCavity}. The corresponding control functions are given by
\[
p(t) = \frac{A}{2} \left(1 + \tanh (\delta ( t - \tau)) \right), \ \ q(t) = 0.
\]
The values of $A, \delta$, and $\tau$ are chosen as  \SI{10}{\mega\hertz}, \SI{-0.05}{\giga\hertz} and \SI{2000}{\nano\second} for the qudit and \SI{15}{\mega\hertz}, \SI{-0.1}{\giga\hertz} and \SI{200}{\nano\second} for the resonator. We simulate the problem up to the final time \(T=\SI{2500}{\nano\second}\). We first study the convergence of the proposed schemes, then demonstrate the low-rank method, and finally present an example showing that trace renormalization, while preventing blow-up, may alter the accuracy of the solution. To obtain the reference solution,  we use a large number of time steps for fourth-order scheme with implicit flow operator. 

\begin{table}[]
\centering 
\caption{Parameters for the qudit-cavity problem.}
\label{table:parameters_QuditCavity}
\begin{tabular}{|c|c|c|c|c|}
\hline  
$\omega_0/2\pi$ \SI{}{\giga\hertz} & $\omega_1/2\pi$ \SI{}{\giga\hertz} & $\xi_0/2\pi$ \SI{}{\giga\hertz} & $\xi_1/2\pi$ \SI{}{\mega\hertz} & $J_{01}$ \SI{}{\giga\hertz} \\
\hline  
$4.41666$ & $6.84081$ &  $0.23056$& $0$ &0 \\
\hline  
\hline  
$\xi_{01}$ \SI{}{\mega\hertz} & $T_{1}^0$ \SI{}{\micro\second} &$T_{1}^1$ \SI{}{\micro\second} & $T_{2}^0$ \SI{}{\micro\second} &$T_{2}^1 $ \SI{}{\micro\second}\\
\hline  
 $1.176$&80&0.3892&26& $\infty$  \\
\hline  
\end{tabular}
\end{table}

\begin{figure}[htb]
\centering
\includegraphics[width=1.0\textwidth, trim=0mm 6mm 0mm 0mm, clip=true]
{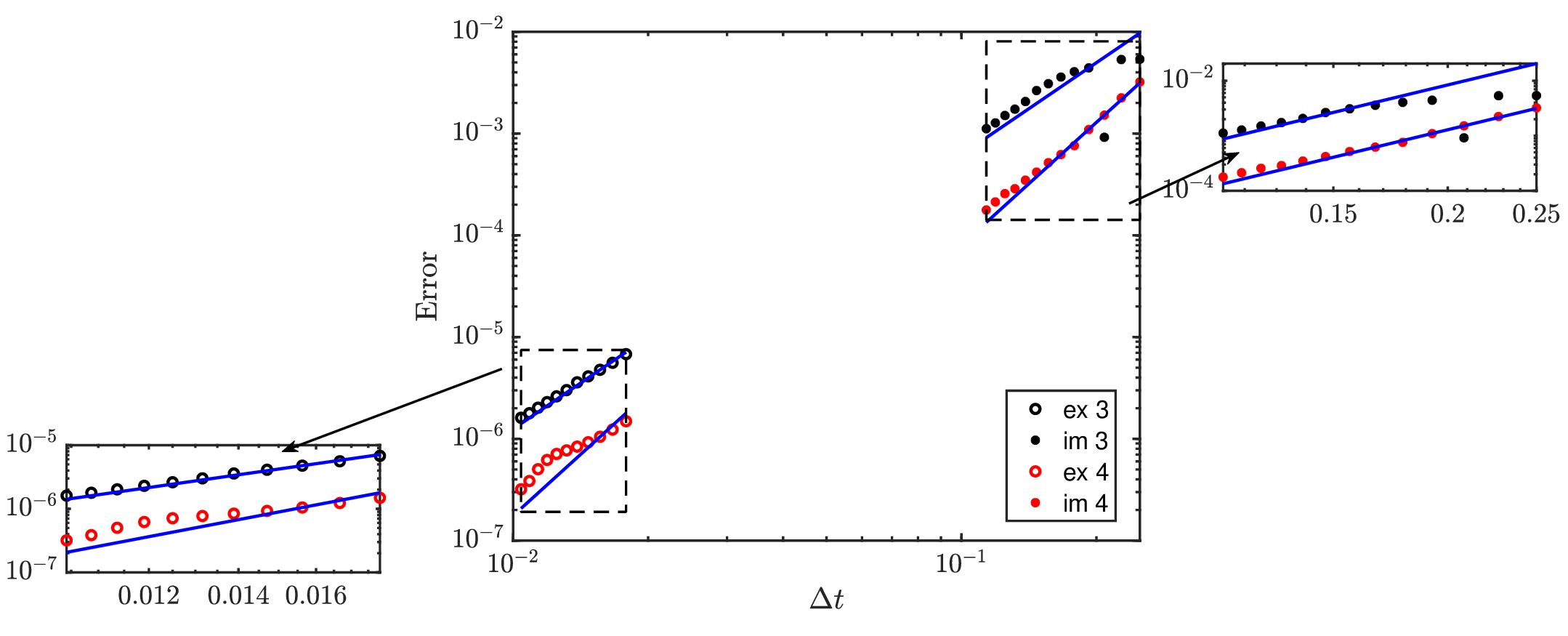}
\caption{Errors as functions of the time step size for the third- and fourth-order schemes with explicit and implicit flow operators in the qudit-resonator problem.}
\label{fig:error_QuditCavity}
\end{figure}

To study the convergence behavior, we consider the third- and fourth-order schemes with both explicit and implicit flow operators. The error is measured in the Frobenius norm of the density matrix at the final time, and Figure~\ref{fig:error_QuditCavity} plots the error versus the time step size. The observed convergence rates are consistent with the theoretical orders, except for the fourth-order scheme with an explicit flow operator, whose rate is slightly below the expected one. This is likely due to the fact that the error is already very small and hence influenced by machine precision. In addition, for a prescribed error level, the schemes with implicit flow operators require substantially fewer time steps than the corresponding explicit schemes. We also note that, for each method, the largest time step size used in the figure is close to the empirical maximum value for which trace renormalization does not noticeably deteriorate the accuracy. Overall, Figure~\ref{fig:error_QuditCavity} shows that the schemes with implicit flow operators offer greater flexibility in the choice of time step size.

\begin{figure}[]
\centering
\includegraphics[width=0.49\textwidth, trim=18mm 68mm 20mm 70mm, clip=true]{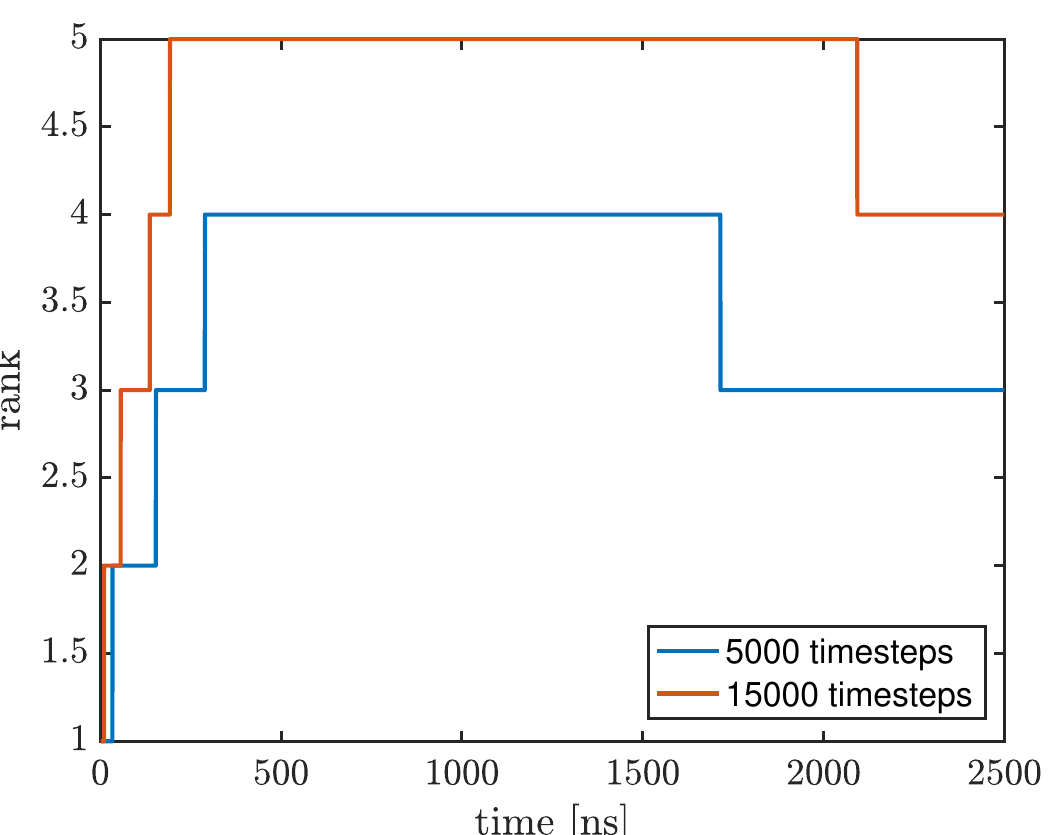}
\includegraphics[width=0.49\textwidth, trim=20mm 68mm 20mm 70mm, clip=true]{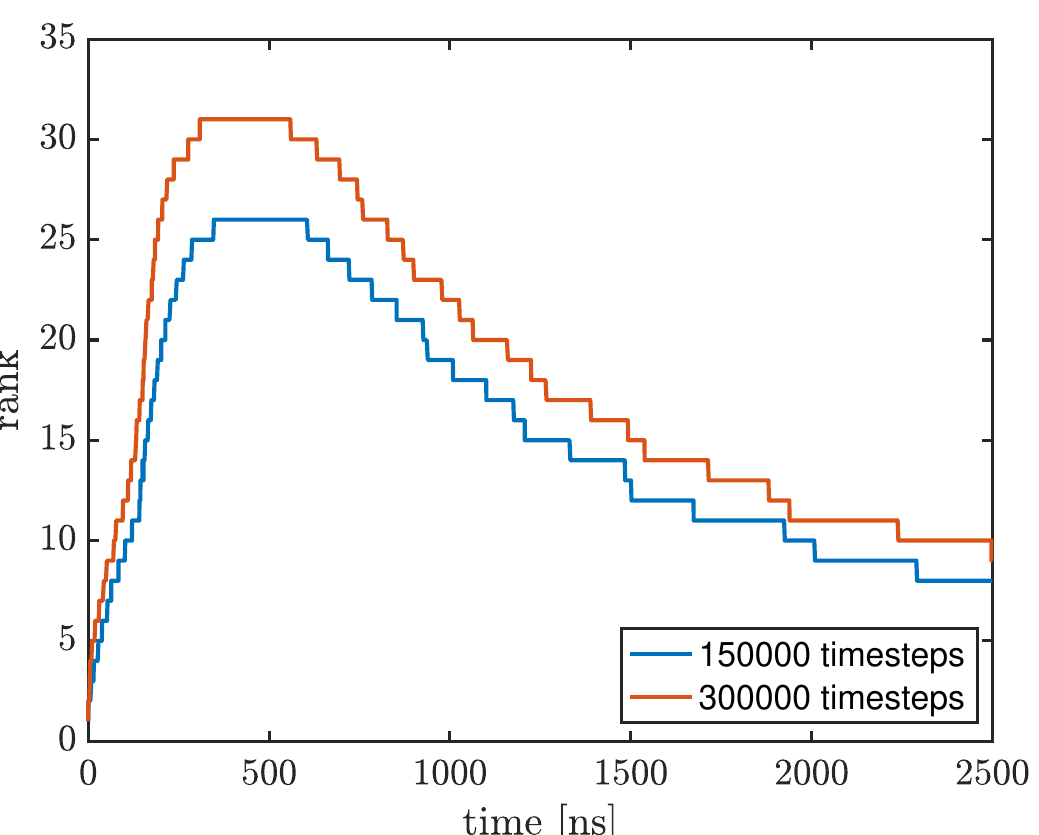}
\caption{Rank evolution for the fourth-order scheme with implicit (left) and explicit (right) flow operators.}
\label{fig:rank_QuditCavity}
\end{figure}

Figure~\ref{fig:rank_QuditCavity} displays the rank history for the fourth-order scheme. Here we take \(\kappa=0.33\) for the implicit-flow schemes and \(\kappa=0.53\) for the explicit-flow schemes. These parameters are chosen empirically so that the truncation error is comparable to the discretization error before truncation. For both types of flow operators, the final ranks remain similar across different time step sizes. While finer time steps may lead to temporarily larger intermediate ranks, truncation reduces them effectively, and the final rank increases only slightly. Together with Figure~\ref{fig:error_QuditCavity}, this indicates that the scheme with the implicit flow operator is more efficient.

\begin{figure}[]
\centering
\includegraphics[width=0.49\textwidth, trim=9mm 66mm 22mm 74mm, clip=true]{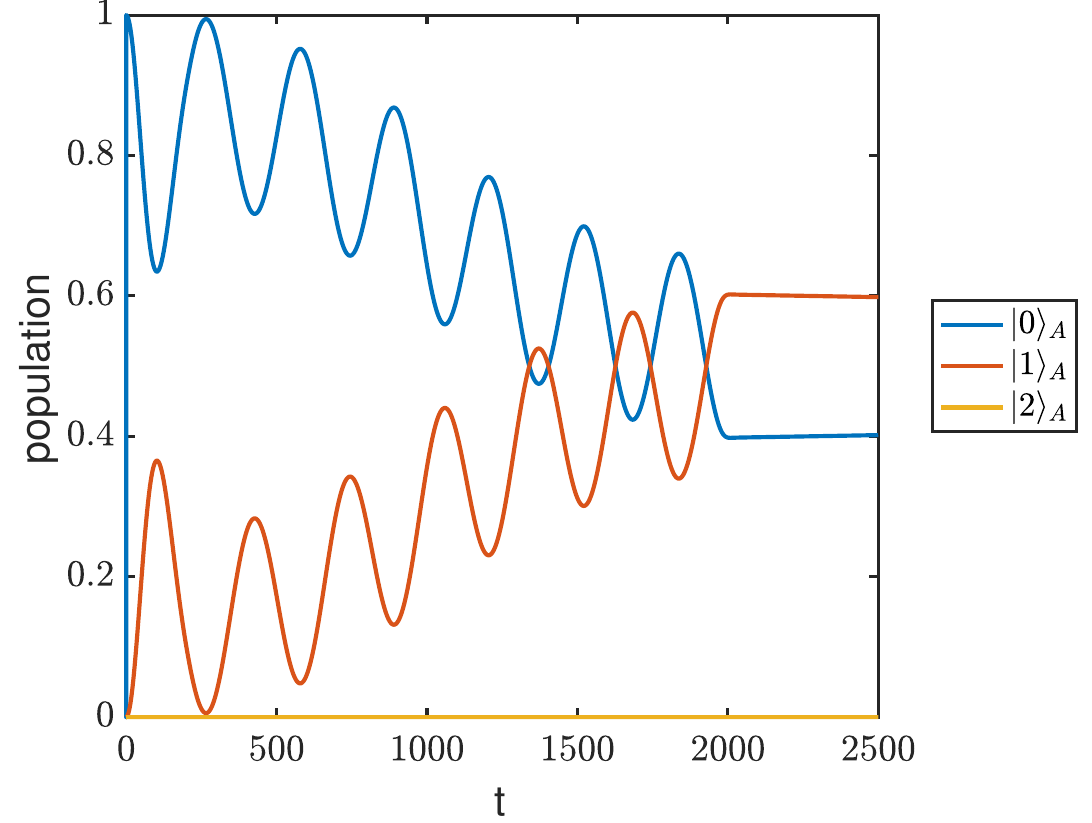}
\includegraphics[width=0.49\textwidth, trim=9mm 66mm 22mm 74mm, clip=true]{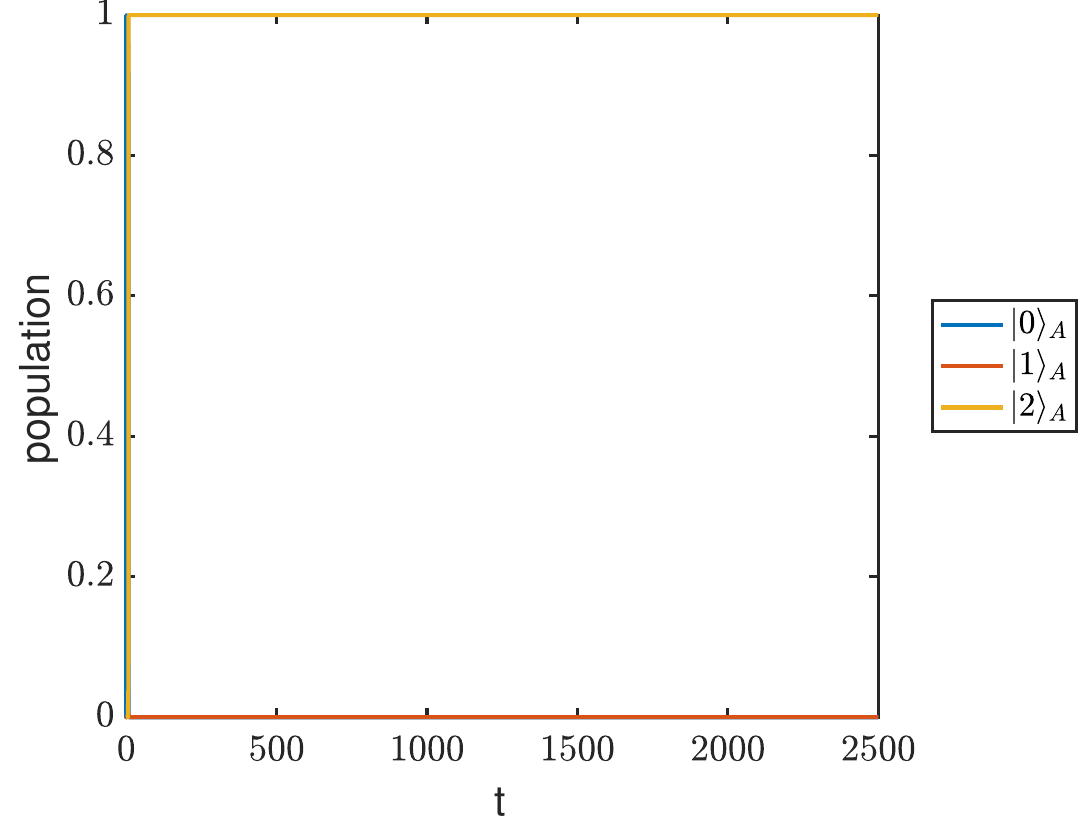}
\caption{Population in the qudit: reference solution (left) and fourth-order explicit scheme with $10000$ time steps (right).}
\label{fig:QuditCavity_population}
\end{figure}

Finally, we investigate the effect of renormalization on the numerical solution. Figure~\ref{fig:QuditCavity_population} compares the population of the qudit obtained with a very large number of time steps (left), which serves as a reference solution, and that produced by the fourth-order explicit scheme with $10000$ time steps (right). The right figure shows that renormalization keeps the solution stable, but does not guarantee sufficient accuracy. This is one of the main reasons for studying the stability of the scheme without renormalization in Section~\ref{sec: stability}.

\subsection{Two-Qudit-Resonator System with a CNOT Gate}

\begin{figure}[]
\centering
\includegraphics[width=0.99\textwidth, trim=20mm 1mm 6mm 1mm, clip=true]{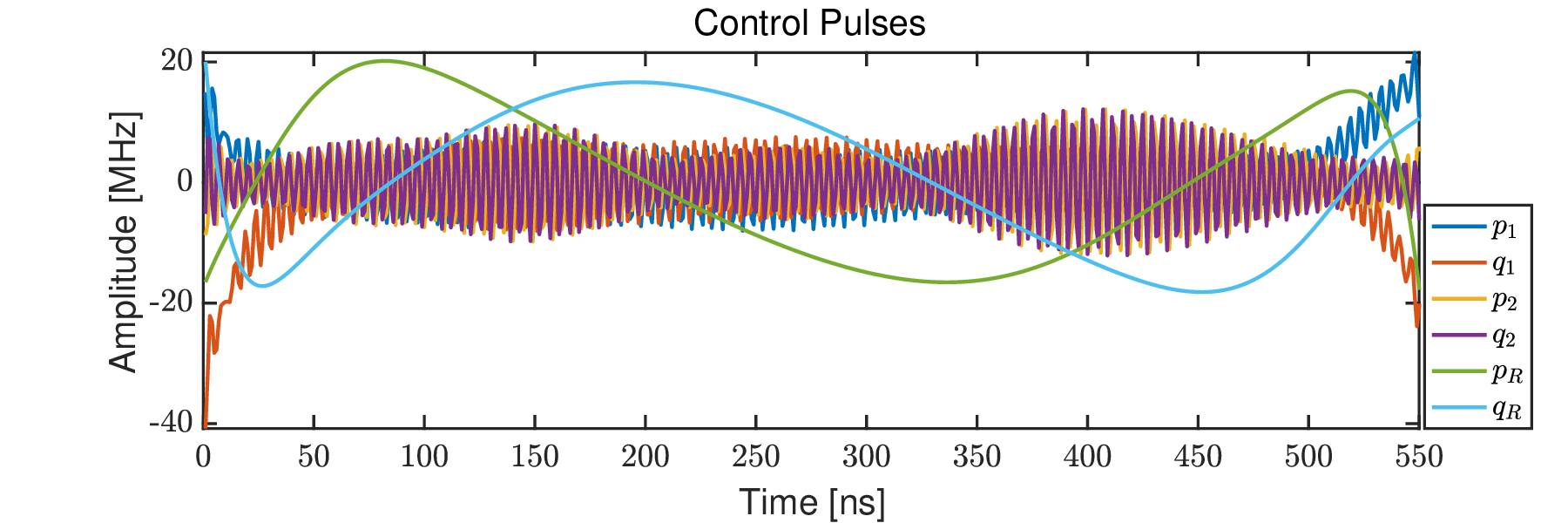}
\caption{Control Pulses for CNOT gate. \label{fig:control_func}}
\end{figure}

In this experiment, we consider a composite quantum system consisting of two qudits, each with four energy levels, and a resonator with ten energy levels, in both closed- and open-system settings. We adopt the same physical parameters and the same CNOT control pulse as in \cite{lee2026high}. The parameter values are listed in Table~\ref{table:parameters_2qudit_resonator}, and the control pulse is shown in Figure~\ref{fig:control_func}. In the closed-system case, the dynamics are governed by the von Neumann equation, that is, without the Lindbladian term. Our primary goal is to assess whether the control pulse from \cite{lee2026high} remains effective for the von Neumann equation and whether its performance is preserved when higher energy levels are included. We also examine how this pulse, originally designed for a closed system, performs in the corresponding open-system setting.

\begin{table}[H]
\centering
\caption{Parameters for the two-qudit--resonator problem. All values are given in \si{\giga\hertz}.}
\label{table:parameters_2qudit_resonator}
\begin{tabular}{|c|c|c|c|c|}
\hline
$\omega_1/2\pi$ & $\omega_2/2\pi$ & $\omega_R/2\pi$ & $\xi_1/2\pi$ & $\xi_2/2\pi$\\
\hline
$4.11$ & $4.82$ & $7.84$ & $2.20(-1)$ & $2.25(-1)$\\
\hline\hline
$\xi_R/2\pi$ & $\xi_{12}/2\pi$ & $\xi_{1R}/2\pi$ & \multicolumn{2}{c|}{$\xi_{2R}/2\pi$}\\
\hline
$2.83(-5)$ & $1.00(-6)$ & $2.49(-3)$ & \multicolumn{2}{c|}{$2.52(-3)$}\\
\hline
\end{tabular}
\end{table}

\begin{figure}[]
\centering
\includegraphics[width=0.49\textwidth, trim=13mm 66mm 27mm 66mm, clip=true]{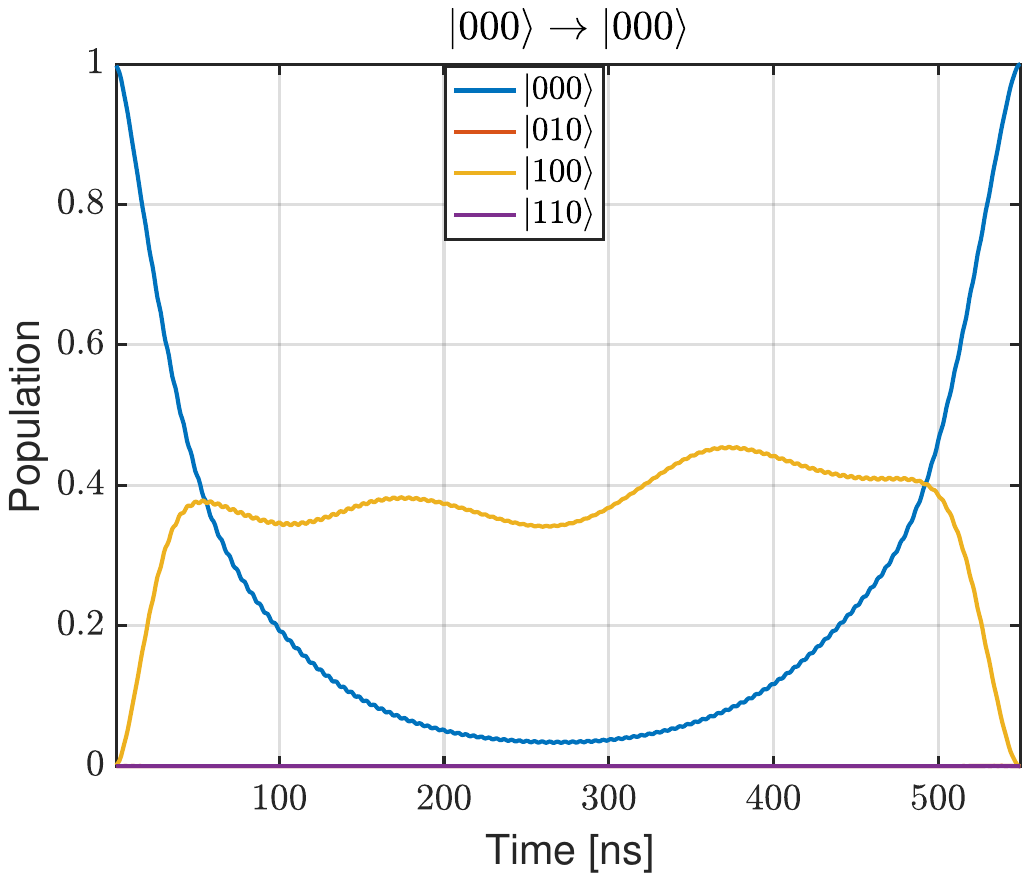}
\includegraphics[width=0.49\textwidth, trim=13mm 66mm 27mm 66mm, clip=true]{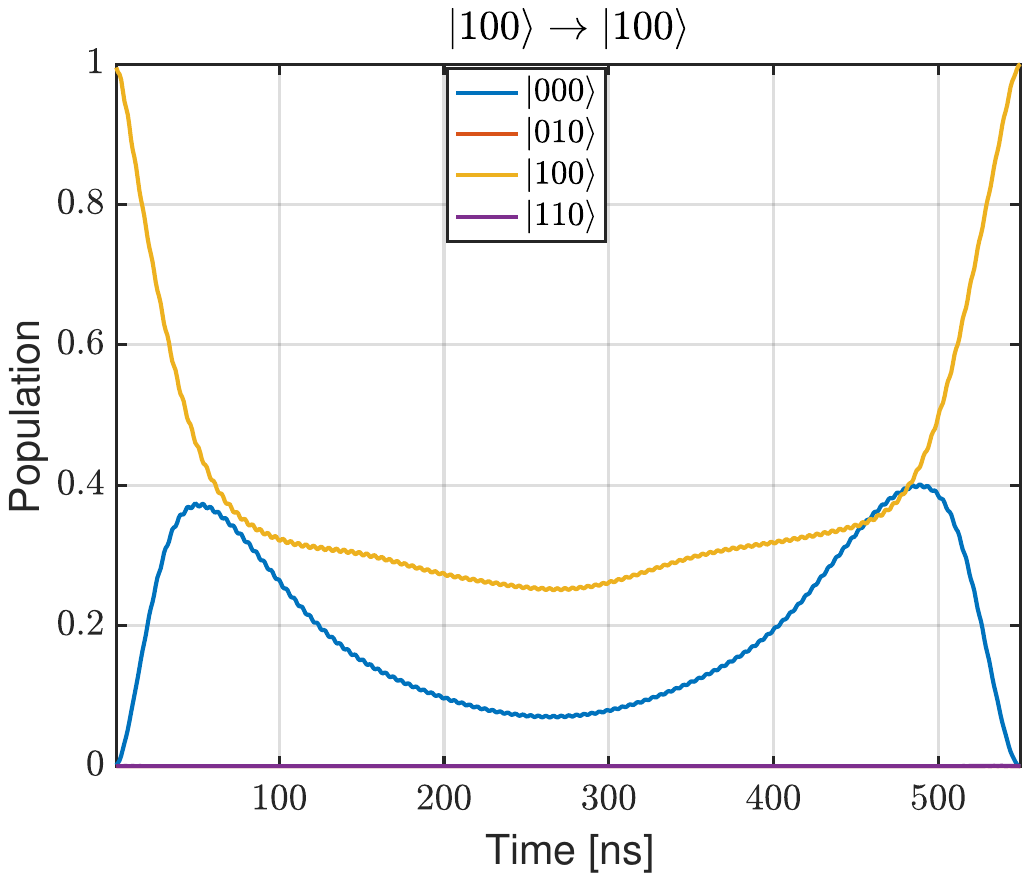}
\includegraphics[width=0.49\textwidth, trim=13mm 66mm 27mm 66mm, clip=true]{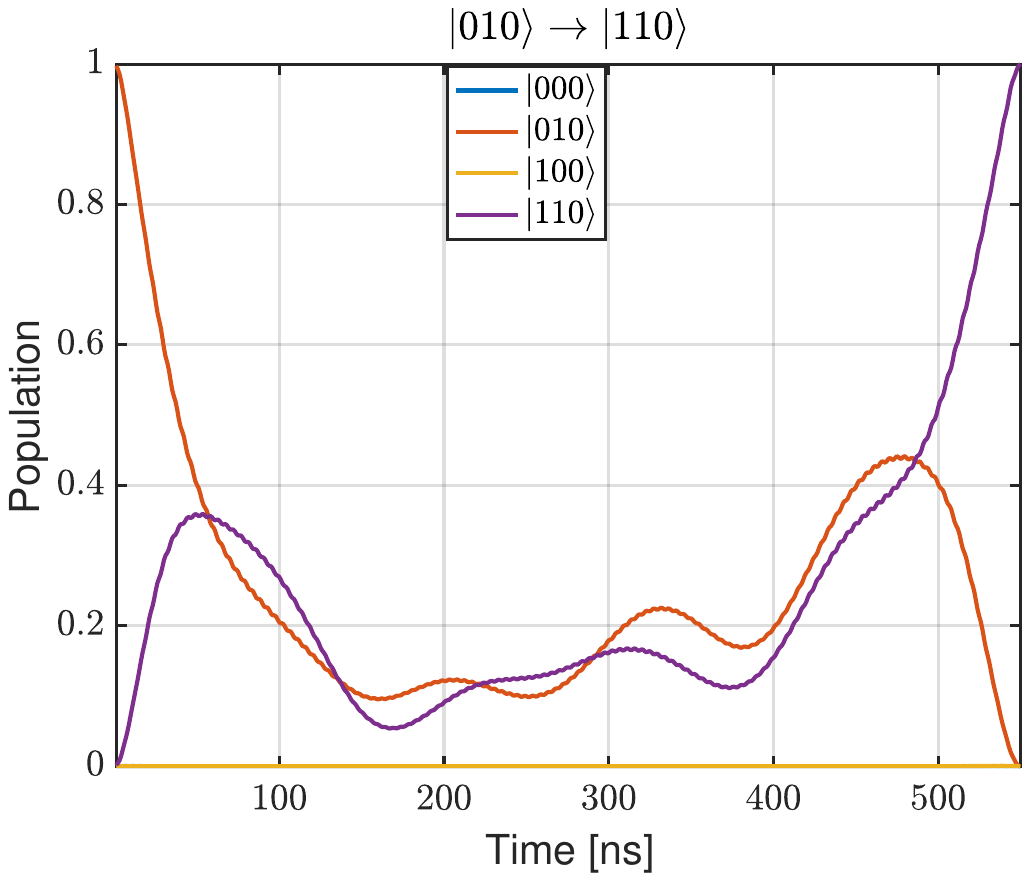}
\includegraphics[width=0.49\textwidth, trim=13mm 66mm 27mm 66mm, clip=true]{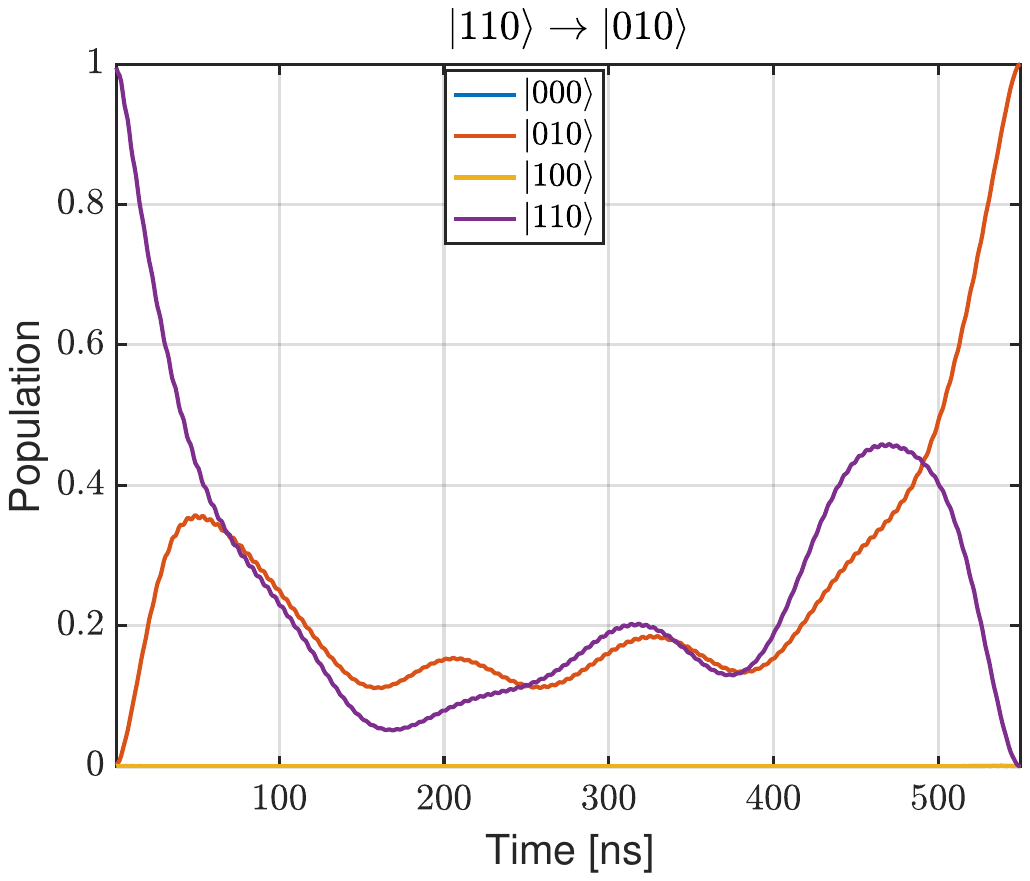}
\caption{Populations with pure initial states in closed system. \label{fig:Population_2QuditsResonator_pure_state}}
\end{figure}

We first examine the performance of the CNOT gate in the closed-system setting governed by the von Neumann equation. We consider four pure initial states,
\[
\ket{000},\quad \ket{010},\quad \ket{100},\quad \ket{110}.
\]
Figure~\ref{fig:Population_2QuditsResonator_pure_state} shows that the gate flips qudit 1 between the ground and excited states when qudit 2 is in the excited state, while leaving qudit 1 unchanged when qudit 2 is in the ground state. This behavior is consistent with the intended action of the CNOT gate.

We next increase the number of resonator energy levels from 10 to 20 while keeping each qudit at four levels. We also consider a larger configuration with 6 energy levels for each qudit and 20 energy levels for the resonator. In all cases, the control pulse exhibits the same qualitative behavior. For brevity, we omit the corresponding figures, as they are qualitatively similar to Figure~\ref{fig:Population_2QuditsResonator_pure_state}.

\begin{figure}[]
\centering
\includegraphics[width=0.49\textwidth, trim=13mm 66mm 27mm 70mm, clip=true]{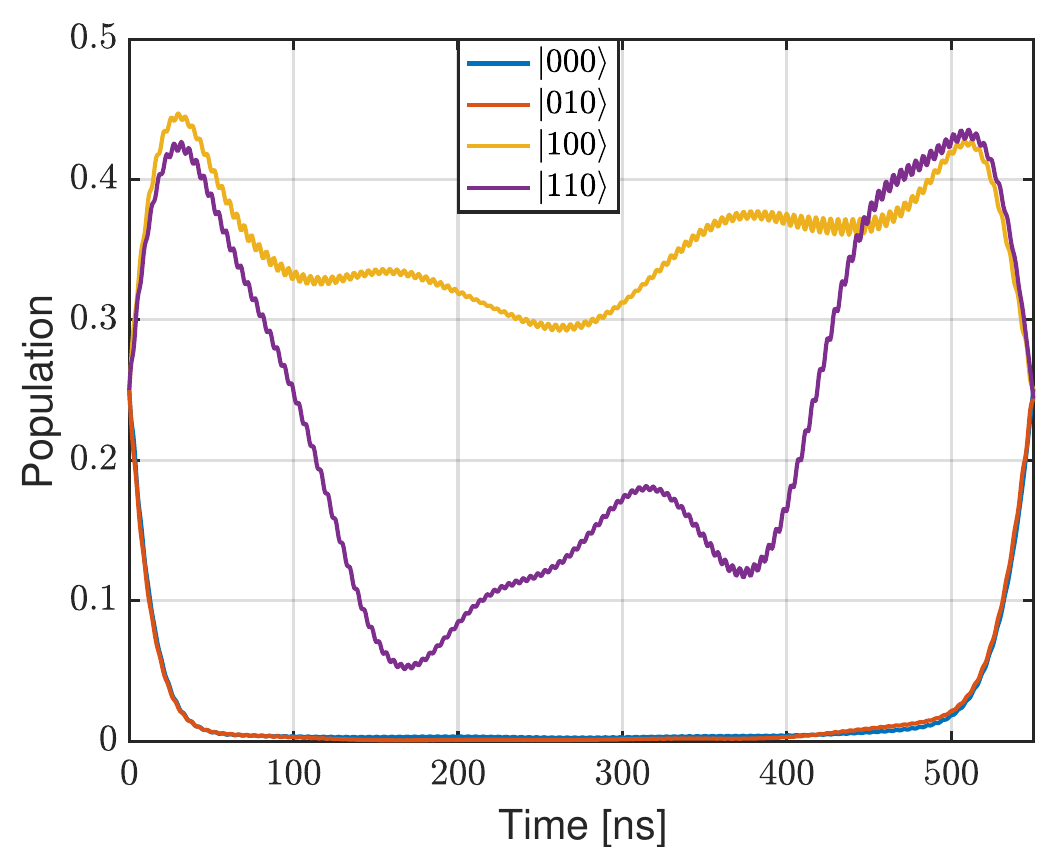}
\caption{Populations  with pure  state in open system.\label{fig:Population_2QuditsResonator_Lindblad_pure_state}}
\end{figure}

Next, we test the same CNOT control pulse in the open-system setting. The relaxation and dephasing times for the two qudits and the resonator are chosen as
\[
T_{1}^{1}=T_{1}^{2}=40~\si{\micro\second},\qquad
T_{1}^{R}=200~\si{\micro\second},\qquad
T_{2}^{1}=T_{2}^{2}=20~\si{\micro\second},\qquad
T_{2}^{R}=500~\si{\micro\second}.
\]
These values are experimentally plausible. The initial density matrix is taken to be
\[
\rho_0=\ket{\psi}\bra{\psi},
\qquad
\ket{\psi}
=
\frac{1}{2}\bigl(\ket{000}+\ket{010}+\ket{100}+\ket{110}\bigr).
\]
We compute the solution using the fourth-order scheme with implicit flow operator and \(\Delta t=0.01\). The corresponding population dynamics are shown in Figure~\ref{fig:Population_2QuditsResonator_Lindblad_pure_state}. Let \(P_{ijk}(t)\) denote the population of the state \(\ket{ijk}\) at time \(t\). At the final time, the numerical results give
\[
|P_{000}(550)-P_{000}(0)| = 2.79\times 10^{-4},\qquad
|P_{010}(550)-P_{110}(0)| = 5.57\times 10^{-3},
\]
\[
|P_{110}(550)-P_{010}(0)| = 5.54\times 10^{-3},\qquad
|P_{100}(550)-P_{100}(0)| = 1.37\times 10^{-4}.
\]

These results indicate that, for the above parameter regime, the control pulse remains effective in realizing the CNOT gate in the open-system setting.

\section{Conclusion}\label{sec:conclusion}
In this work, we developed novel high-order low-rank structure-preserving numerical schemes for the Lindblad equation. The proposed methods build on NPI framework with two key ideas. First, Gaussian quadrature is replaced by Gregory quadrature, which allows the use of equidistant nodes while maintaining high-order accuracy. Second, the approximation of the flow operator is fixed at the highest available order within the construction, which further improves efficiency in practical computations. As a result, the new schemes retain the main advantages of the previous framework, including high-order accuracy, low-rank structure, and preservation of the key physical properties of the density matrix, while reducing the overall computational cost.

We also investigated the stability of the resulting schemes without trace renormalization. For a representative single-qubit model problem, the analysis clarifies the role of the flow operator approximation and shows how the stability behavior depends on the dissipative parameters and the time step size. In particular, for sufficiently large dephasing times, the stability region is determined primarily by the stability of the flow operator. This provides useful guidance for the practical implementation of the proposed methods.
Low-rank truncation is incorporated into the  framework. The numerical experiments demonstrate that the proposed methods achieve the expected convergence orders and perform well on several representative problems, including a two-qubit problem with known exact solution, a controlled qudit-resonator system, and CNOT gate simulations for two qudits coupled with a resonator in both closed- and open-system settings. 

\appendix
\section*{Appendix}
In the Appendix, we establish the proof of the stability results for the fourth-order scheme. The analysis is carried out by first vectorizing the matrices and then reformulating the multistep scheme as an equivalent one-step method. The stability is then characterized by the amplification matrix of the resulting one-step system. 

\section{Proof of Theorem \ref{thm: order 4}}\label{sec:stability_proof_order_4}
By vectorizing the matrices $\rho^{n+5,k}$, we can rewrite the fourth-order approximation in \eqref{eqn: approximation_order_4} as follows:

\begin{align*}
\boldsymbol{\rho}^{n+5,k}
&= P_5^{(4)}\boldsymbol{\rho}^n
+ \Delta t\Big(
\omega^{(4)}_0 Q_{U^{(4)}L,5}\boldsymbol{\rho}^n
+ \omega^{(4)}_1 Q_{U^{(4)}L,4}\boldsymbol{\rho}^{n+1}
+ \omega^{(4)}_2 Q_{U^{(4)}L,3}\boldsymbol{\rho}^{n+2} \\
&\quad
+ \omega^{(4)}_3 Q_{U^{(4)}L,2}\boldsymbol{\rho}^{n+3}
+ \omega^{(4)}_4 Q_{U^{(4)}L,1}\boldsymbol{\rho}^{n+4}
\Big)
+ \Delta t\,\omega^{(4)}_5 Q_L \boldsymbol{\rho}^{n+5,k-1},
\end{align*}
for $k=1,2,3,4$. Introducing the augmented vector
\[
\bm{\boldsymbol{\rho}^{n+5}\\ \boldsymbol{\rho}^{n+4}\\ \boldsymbol{\rho}^{n+3}\\ \boldsymbol{\rho}^{n+2}\\ \boldsymbol{\rho}^{n+1}}
:=
\bm{\boldsymbol{\rho}^{n+5,4}\\ \boldsymbol{\rho}^{n+4}\\ \boldsymbol{\rho}^{n+3}\\ \boldsymbol{\rho}^{n+2}\\ \boldsymbol{\rho}^{n+1}},
\]
the scheme can be formulated  as the one-step system
\[
 \bm{\boldsymbol{\rho}^{n+5}\\\boldsymbol{\rho}^{n+4}\\\boldsymbol{\rho}^{n+3}\\\boldsymbol{\rho}^{n+2}\\\boldsymbol{\rho}^{n+1}}
=
 \bm{A_1&A_2&A_3&A_4&A_5\\
I&O&O&O&O\\
O&I&O&O&O\\
O&O&I&O&O\\
O&O&O&I&O}
 \bm{\boldsymbol{\rho}^{n+4}\\\boldsymbol{\rho}^{n+3}\\\boldsymbol{\rho}^{n+2}\\\boldsymbol{\rho}^{n+1}\\\boldsymbol{\rho}^{n}},
\]
where
\begin{align*}
    A_1 &=\left(I+\omega^{(4)}_5\Delta tQ_L+(\omega^{(4)}_5\Delta tQ_L)^2+(\omega^{(4)}_5\Delta tQ_L)^3\right) \omega^{(4)}_4 \Delta tQ_{U^{(4)}L,1}+ (\omega^{(4)}_5\Delta tQ_L)^4,\\
    A_k &= \left(I+\omega^{(4)}_5\Delta tQ_L+(\omega^{(4)}_5\Delta tQ_L)^2 +(\omega^{(4)}_5\Delta tQ_L)^3\right) \omega^{(4)}_{5-k} \Delta tQ_{U^{(4)}L,k}, \quad \text{for }k = 2,3,4,\\  
    A_5 &=\left(I+\omega^{(4)}_5\Delta tQ_L+(\omega^{(4)}_5\Delta tQ_L)^2
    +(\omega^{(4)}_5\Delta tQ_L)^3\right) \left( P^{(4)}_5 +\omega^{(4)}_0 \Delta tQ_{U^{(4)}L,5} \right) .
\end{align*}
 
Define
\begin{align*}
    C&  := \omega^{(4)}_5 \frac{\Delta t}{T_1} \sum_{j=0}^2 \left(\omega^{(4)}_5\right)^j\left(\frac{\Delta t}{T_2}\right)^j, \\
    B^{(4)}& := \sum_{j=0}^3 \left(\omega^{(4)}_5\right)^j\left(\frac{\Delta t}{T_2}\right)^j.
\end{align*}

Then we have
\begin{align*}
    A_1 &= \omega^{(4)}_4\bm{0&0&0& \frac{\Delta t}{T_1} +\frac{\Delta t}{T_2} 
    |\varsigma^{(4)} ( 1)|^2 C +\cfrac{(\omega^{(4)}_5\Delta t )^4} {T_1(T_2)^3}
    \\0&0&0&0
    \\0&0&0&0
    \\0&0&0&  \frac{\Delta t}{T_2}|\varsigma^{(4)} ( 1)|^2
    B^{(4)}+\cfrac{(\omega^{(4)}_5\Delta t )^4} {(T_2)^4}}
    \\
    A_k &=  \omega^{(4)}_{5-k} \bm{0&0&0&\frac{\Delta t}{T_1} +\frac{\Delta t}{T_2} 
    |\varsigma^{(4)} ( k)|^2 C\\0&0&0&0\\0&0&0&0\\0&0&0&\frac{\Delta t}{T_2}|\varsigma^{(4)} ( k)|^2
    B^{(4)}}, \text{    for  } k=2,3,4\\
    A_5 &
    =\bm{1&0&0&\overline{\varsigma^{(4)} ( 5) }C +\omega^{(4)}_0 \frac{\Delta t}{T_1} +\omega^{(4)}_0\frac{\Delta t}{T_2} 
    |\varsigma^{(4)} ( 5)|^2 C
    \\0&\overline{\varsigma^{(4)} ( 5)}&0&0
    \\0&0&\varsigma^{(4)} ( 5)&0\\0&0&0&|\varsigma^{(4)} ( 5)|^2B^{(4)} +\omega^{(4)}_0\frac{\Delta t}{T_2}|\varsigma^{(4)} ( 5)|^2
    B^{(4)}}.
\end{align*}
The matrix $$\bm{A_1&A_2&A_3&A_4&A_5\\
I&O&O&O&O\\
O&I&O&O&O\\
O&O&I&O&O\\
O&O&O&I&O},$$ 
is similar to the matrix
$$\bm{ 
 O   & I   & O   & O   & O\\
 O   & O   & I   & O   & O\\
 O   & O   & O   & I   & O\\
 O   & O   & O   & O   & I\\
 A_5 & A_4 & A_3 & A_2 & A_1}.$$
 Consequently, the two matrices share the same eigenvalues. The latter matrix is the block companion matrix of the matrix polynomial
\[
L(\lambda)= \lambda^5 - A_1\lambda^4 - A_2\lambda^3 - A_3\lambda^2 - A_4\lambda - A_5 .
\]

Since each of the five matrices \( A_1, \ldots, A_5\) is upper triangular, the matrix polynomial \(L(\lambda)\) is also upper triangular. Consequently, the diagonal entries of \(L(\lambda)\) are
\[
\lambda^5-1,\qquad \lambda^5-\overline{\varsigma^{(4)}(5)},\qquad \lambda^5-\varsigma^{(4)}(5),
\]
and
\[
\lambda^5- \sum_{k=1}^{5}\lambda^{5-k}\omega^{(4)}_{5-k}  \frac{\Delta t}{T_2}|\varsigma^{(4)}(k)|^2 B^{(4)}
-\lambda^4 \left(\cfrac{\omega^{(4)}_5\Delta t}{T_2}\right)^4
-|\varsigma^{(4)}(5)|^2 B^{(4)}.
\]
Therefore, the stability condition is characterized by the roots of these diagonal entries, which proves Theorem \ref{thm: order 4}.

\section{Proof of Theorem \ref{thm: order p}}\label{sec:stability_proof_order_p}
After vectorizing the matrices \(\rho^{n+n_p,k}\), we express the $p$-th-order approximation in \eqref{eq:NPI_general_new} in the following form:
\begin{align*}
\boldsymbol{\rho}^{n+n_p,0} & = \boldsymbol{\rho}^{n+n_p-1},  \\
\boldsymbol{\rho}^{n+n_p,k} &=  P_{n_p}^{(p)}\boldsymbol{\rho}^n + \Delta t \sum_{i=0}^{2p-4} \omega^{(p)}_iQ_{U^{(p)}L,n_p-i}\boldsymbol{\rho}^{n+i} + \Delta t\omega^{(p)}_{n_p} Q_L \boldsymbol{\rho}^{n+n_p,k-1},
\end{align*}
for $k=1,2,\cdots,p$.  Introducing the augmented vector
\[
\bm{\boldsymbol{\rho}^{n+n_p}\\\boldsymbol{\rho}^{n+n_p-1} \\\vdots \\\boldsymbol{\rho}^{n+2}\\\boldsymbol{\rho}^{n+1} }
:=
\bm{\boldsymbol{\rho}^{n+n_p,p}\\\boldsymbol{\rho}^{n+n_p-1} \\\vdots \\\boldsymbol{\rho}^{n+2}\\\boldsymbol{\rho}^{n+1} },
\]
the scheme can be rewritten  as the one-step system

 \begin{align*}
     \bm{\boldsymbol{\rho}^{n+n_p}\\\boldsymbol{\rho}^{n+n_p-1} \\\vdots \\\boldsymbol{\rho}^{n+3}\\\boldsymbol{\rho}^{n+2}\\\boldsymbol{\rho}^{n+1} }
     &= \bm{A_1&A_2&\cdots&A_{n_p-2}&A_{n_p-1}&A_{n_p}\\
I&O&\cdots&O&O&O\\
O&I&\cdots&O&O&O\\
\vdots&\vdots&\ddots &\vdots&\vdots&\vdots\\
O&O&\cdots&I&O&O\\
O&O&\cdots&O&I&O}\bm{\boldsymbol{\rho}^{n+n_p-1} \\\boldsymbol{\rho}^{n+n_p-2}\\\vdots \\\boldsymbol{\rho}^{n+2} \\\boldsymbol{\rho}^{n+1} \\\boldsymbol{\rho}^{n}},
 \end{align*}
  
where
\begin{align*}
    A_1 &=\left(I+\sum_{i=1}^{p-1} (\omega^{(p)}_{n_p}\Delta tQ_L)^i
    \right) \omega^{(p)}_{n_p-1} \Delta tQ_{U^{(p)}L,1}
    + \left(\omega^{(p)}_{n_p}\Delta tQ_L\right)^p, \\
    A_k &= \left(I+\sum_{i=1}^{p-1} (\omega^{(p)}_{n_p}\Delta tQ_L)^k
    \right) \omega^{(p)}_{n_p-k} \Delta tQ_{U^{(p)}L,k}, \text{  for  } k = 2,3,\cdots,n_p-1,\\
    A_{n_p} &= \left(I+\sum_{k=1}^{p-1} (\omega^{(p)}_{n_p}\Delta tQ_L)^k
    \right) (P^{(p)}_{n_p}+ \omega^{(p)}_0 \Delta tQ_{U^{(p)}L,n_p}).
\end{align*}
Define
\begin{align*}
    C^{(p)} : & = \sum_{k=0}^{p-2} \dfrac{\omega^{(p)}_{n_p}\Delta t}{T_1}\left(\dfrac{\omega^{(p)}_{n_p}\Delta t}{T_2} \right)^{k-1},  \\
     B^{(p)}: & = \sum_{k=0}^{p-1} \left(\dfrac{\omega^{(p)}_{n_p}\Delta t}{T_2} \right)^k.
\end{align*}
Then
\begin{align*}
    A_1 &=\omega^{(p)}_{n_p-1} \bm{0&0&0& \frac{\Delta t}{T_1} +\frac{\Delta t}{T_2} 
    |\varsigma^{(p)} ( 1)|^2 C^{(p)} +\cfrac{(\omega^{(p)}_{n_p}\Delta t )^p} {T_1(T_2)^{p-1}}
    \\
    0&0&0&0\\0&0&0&0\\0&0&0& \frac{\Delta t}{T_2}|\varsigma^{(p)} ( 1)|^2
    B^{(p)}+\cfrac{(\omega^{(p)}_{n_p}\Delta t )^p} {(T_2)^p}},
 \\
   A_k &=\omega^{(p)}_{n_p-k}  \bm{0&0&0&\frac{\Delta t}{T_1} +\frac{\Delta t}{T_2} 
    |\varsigma^{(p)} ( k)|^2 C^{(p)}\\
    0&0&0&0\\0&0&0&0\\0&0&0&\frac{\Delta t}{T_2}|\varsigma^{(p)} ( k)|^2
    B^{(p)}}, \text{  for  } k = 2,\cdots,n_p-1,\\
    A_{n_p} &= \bm{1&0&0&\overline{\varsigma^{(p)} (n_p)}C^{(p)} +\omega^{(p)}_0 \frac{\Delta t}{T_1} +\omega^{(p)}_0\frac{\Delta t}{T_2} 
    |\varsigma^{(p)} ( n_p)|^2 C^{(p)}\\0&\overline{\varsigma^{(p)} (n_p)}&0&0\\0&0&\varsigma^{(p)} (n_p)&0\\0&0&0&|\varsigma^{(p)} (n_p)|^2B^{(p)}+\omega^{(p)}_0\frac{\Delta t}{T_2}|\varsigma^{(p)} ( n_p)|^2
    B^{(p)}}
    .
\end{align*}

Let $\mathcal{M}^{(p)}= \bm{A_1&A_2&\cdots&A_{n_p-2}&A_{n_p-1}&A_{n_p}\\
I&O&\cdots&O&O&O\\
O&I&\cdots&O&O&O\\
\vdots&\vdots&\ddots &\vdots&\vdots&\vdots\\
O&O&\cdots&I&O&O\\
O&O&\cdots&O&I&O}$ be the amplification matrix. By a block permutation, $\mathcal{M}^{(p)}$ is equivalent to 
\[
\bm{
 O   &I   & O & \cdots&O&O\\
 O   & O   &I   & \cdots   &O&O\\
 \vdots   &  \vdots   &  \vdots   & \ddots     & O&O
 \\
 O   & O   & O   &\cdots   &I  & O \\
 O   & O   & O   &\cdots     & O   & I\\
 A_{n_p} & A_{n_p-1} & A_{n_p-2} &\cdots  & A_2 & A_1} ,\]
which is the block companion matrix associated with the matrix polynomial
 \[L(\lambda) =  \lambda^{n_p} - \sum_{k=1}^{n_p} A_k\lambda^{n_p-k}.\]
 Since each \(A_k\), \(k=1,2,\dots,n_p\), is upper triangular, the matrix polynomial \(L(\lambda)\) is also upper triangular. Its diagonal entries are
\[
\lambda^{n_p}-1,\qquad
\lambda^{n_p}-\overline{\varsigma^{(p)}(n_p)},\qquad
\lambda^{n_p}-\varsigma^{(p)}(n_p),
\]
and
\[
\lambda^{n_p}
-\sum_{k=1}^{n_p} \lambda^{n_p-k}\Delta t\,\omega^{(p)}_{n_p-k}\frac{\Delta t}{T_2}|\varsigma^{(p)} ( k)|^2
    B^{(p)}
-\lambda^{n_p-1}\bigl(\omega^{(p)}_{n_p}\Delta t\bigr)^p\frac{1}{(T_2)^p}
-\left| \varsigma^{(p)}(n_p)\right|^2 B^{(p)}.
\]
Therefore, the stability condition is determined by the roots of these diagonal polynomials. This completes the proof of Theorem \ref{thm: order p}.
%

%
%
\bibliographystyle{siam}
\bibliography{ref}
\end{document}